\documentclass[11pt, a4paper, english]{article}

\usepackage{amsmath, amscd, latexsym, babel, amstext}
\usepackage{amssymb}
\usepackage{enumerate}
\usepackage{theorem}
\usepackage{amsfonts}
\usepackage[a4paper,margin=2cm,vmargin={1cm,2cm},includeheadfoot]{geometry}
\usepackage[scaled=0.92]{helvet}
\usepackage[latin 1]{inputenc}

\newcommand{\R}{{\mathbb {R}}}
\newcommand{\prob}{\mathbb{P}}
\newcommand{\E}{{\mathbb{E}}}

\newcommand{\D}{{\mathbb {D}}}
\newcommand{\Le}{{\mathbb {L}}}
\newcommand{\Q}{{\mathbb {Q}}}

\theoremheaderfont{\bfseries}
\newtheorem{theorem}{Theorem}
\newtheorem{example}[theorem]{Example}
\newtheorem{definition}[theorem]{Definition}
\newtheorem{proposition}[theorem]{Proposition}
\newtheorem{lemma}[theorem]{Lemma}
\newtheorem{corollary}[theorem]{Corollary}
{\theorembodyfont{\rmfamily}
\newtheorem{remark}[theorem]{Remark}}
\newenvironment{proof}[1][Proof]{\textbf{#1.} }{\ \rule{0.5em}{0.5em}}

\numberwithin{equation}{section}
\numberwithin{theorem}{section}
\title{A Malliavin-Skorohod calculus in $L^0$ and $L^1$ \\
for additive and Volterra-type processes}

\author{Giulia Di Nunno\thanks{Department of Mathematics, University of Oslo, P.O. Box 1053 Blindern, 0316 Oslo, and Department of Business and Administration, NHH, Helleveien 30, 5045 Bergen, Norway. Email: giulian@math.uio.no } 
\and Josep Vives\thanks{Facultat de Matem\`atiques, 
Universitat de Barcelona, Gran Via 585, 08007 Barcelona (Catalunya), Spain. E-mail: josep.vives@ub.edu} \thanks{%
Supported by grants MEC MTM 2012 31192 and MEC MTM 2013 40782 P} 
}

\date{\today}
\begin{document}

\maketitle

\begin{abstract}
In this paper we develop a Malliavin-Skorohod type calculus for additive processes in the $L^0$ and 
$L^1$ settings, extending the probabilistic interpretation of the Malliavin-Skorohod operators 
to this context. We prove calculus rules and obtain a generalization of the Clark-Hausmann-Ocone formula for random variables in $L^1$.
Our theory is then applied to extend the stochastic integration with respect to volatility modulated L\'evy-driven Volterra processes recently introduced in the literature. Our work yields to substantially weaker conditions that permit to cover integration with respect to e.g. Volterra processes driven by $\alpha$-stable processes with $\alpha < 2$. The presentation focuses on jump type processes.

\vspace{0.5cm}

\smallskip

%\noindent \textbf{JEL code:}

%, 91B70}.
\end{abstract}

\section{Introduction}

\hspace{0.5cm} Malliavin-Skorohod calculus for square integrable functionals of an additive process is today a well established topic. K. It\^o proved in \cite{I51} the so-called chaos representation property of square integrable functionals of the Brownian motion. A generalized version of this property in terms of a random measure associated to a L\'evy process was proved by the same author in \cite{I56}. Later, a Malliavin-Skorohod calculus for Gaussian processes strongly based on the chaos representation property was developed. 
We refer the reader to \cite{N06} for the Gaussian Malliavin-Skorohod calculus.       

In \cite{NV90} it was proved that an abstract Malliavin-Skorohod calculus could be established on any Hilbert space with Fock space structure. An analogous abstract framework was described in \cite{H86}. Indeed, during the following years, the Malliavin-Skorohod calculus based on the Fock space structure was developped for the standard Poisson process (see \cite{NV90}), for a pure jump L\'evy process or a Poisson random measure (see \cite{BDLOP} and \cite{L04}), for a general L\'evy process (see \cite{DMOP04}, \cite{P08} and \cite{SUV07}), and for additive processes (see \cite{D08} and \cite{Y08}). We refer the reader to \cite{DOP08} for the Malliavin-Skorohod calculus for L\'evy processes.     

In \cite{NV95} a version of the Malliavin-Skorohod calculus for the standard Poisson process was developed 
on the canonical Poisson space introduced by J. Neveu in \cite{N76}. They defined a difference operator and its adjoint and proved that these operators coincided respectively with the gradient and the divergence operators based on the Fock space structure associated to this process. So, this work puts the basis for a Malliavin-Skorohod type calculus beyond $L^2$ in that context. J. Picard, in \cite{P96a} and \cite{P96b}, extended and developped this theory to the more general context of Poisson random measures. Many of these ideas are nicely reviewed in \cite{P09}.

Later, J. L. Sol\'e, F. Utzet and J. Vives introduced a Neveu-type canonical space for the pure jump part of a L\'evy process and defined an increment quotient operator, which turned out to coincide with the gradient operator based on the corresponding Fock space structure, see \cite{SUV05} and \cite{SUV07}. 
On this basis they developed a Malliavin-Skorohod calculus for L\'evy processes beyond the $L^2$ setting in the canonical space, extending the results of \cite{NV95}. More results in this framework were obtained in \cite{ALV08}.

The purpose of this paper is two folded. First we want to set the basis for a Malliavin-Skorohod calculus for general additive processes, which allows to deal with $L^1$ and $L^0$ functionals of the process. We recall that additive processes can be thought of as L\'evy processes without stationary increments, see \cite{S99}. On one hand we extend substantially the theory of \cite{SUV07}, using also ideas from \cite {P96a} and \cite{P96b}. The results also extend the $L^2$ Malliavin-Skorohod calculus developed in \cite{Y08}. Moreover, taking a different perspective, the Skorohod type integral introduced in this paper, defined for the additive processes, extends the 
It\^o integral in $L^1$ (see e.g. \cite{CK}) to the anticipative framework. 

The second goal of the present paper is to discuss explicit stochastic integral representations in the $L^1$ setting. Indeed we prove various rules of calculus and a new version of the Clark-Hausmmann-Ocone (CHO) formula in the $L^1$ setting. This formula extends on the one hand the $L^2$ CHO formulas for L\'evy processes that can be found in \cite{BDLOP}, \cite{DMOP04}, \cite{D08}, \cite{DOP08}, \cite{P09} and \cite{SUV05}. On the other hand the formula extends the pure Brownian CHO formula in the $L^1$ setting obtained by I. Karatzas, J. Li and D. Ocone in \cite{KOL} and also the formula obtained in \cite{P96a} for integrable functionals of the standard Poisson random measure. Moreover, our formula allows to identify the kernels of the martingale representation for additive processes covering many of the cases treated in \cite{C13}, \cite{D07}.

In a summary the original achievement of this paper is to establish and work with techniques proper of canonical spaces to obtain results of stochastic integration in 
$L^0$ and $L^1$ settings. While the statements of the results may not sound surprising as we try to extend the Malliavin-Skorohod integration scheme, the fact that we have substantially enlarged the very set of integrators and integrands opens up for new possible applications. For example, our theory allows to treat stochastic Malliavin-Skorohod integration with respect to $\alpha$-stable processes when $\alpha<2$, in which case there is no second moment available and, in some cases, not the first either. 
We recall that $\alpha$-stable processes are heavy tailed distributions and they appear e.g. in the analysis of financial time series of returns (see e.g. \cite{CT}) and weather linked securities (see e.g. \cite{AZ13}). Moreover, we apply our theory to extend the integral suggested in \cite{BBPV} for Volterra-type processes. Indeed we can treat the case of processes driven by pure jump additive processes in $L^0$ and $L^1$. These models, called volatility modulated Volterra processes (also part of the family of ambit processes), are a flexible class of models used both in turbulence and in energy finance, where the risks may derive from natural phenomena (e.g. wind) with extreme erratic behaviour. In this case the driving noises are characterized by a large tailed distribution, without second moment, see e.g. \cite{BS03}.

The paper is organized as follows. Section 2 is devoted to preliminaries about additive processes. In Section 3 we present fundamental elements of the $L^2$ Malliavin-Skorohod calculus for additive processes as a point of departure of our work. In Section 4, we extend the canonical space for L\'evy processes developed in \cite{SUV07} to the context of additive processes. In Section 5 we introduce a Malliavin calculus in the $L^0$ and $L^1$ settings for Poisson random measures and in particular for pure jump additive processes. We work in the canonical space and we exploit its structure. Section 6 is dedicated to the CHO formula. Our work focuses on the pure jump case. For what the Brownian component is concerned, we recall the results of \cite{KOL} about the CHO formula in the $L^1$. The integration with respect to pure jump volatility modulated Volterra processes is discussed in Section 7.

\section{Preliminaries about additive processes}

\hspace{0.5cm} Consider a real additive process $X=\{X_t, t\geq 0\}$ defined on a complete probability space $(\Omega, {\cal F}, {\prob}).$ Denote by $\E$ the expectation with respect to $\prob.$ Denote by $\{{\cal F}^X_t, t\geq 0\}$ the completed natural filtration of $X$ and define $ {\cal F}^X :=\vee_{t\geq 0} {\cal F}^X_t$. Recall that an additive process is a process with independent increments, stocastically continuous, null at the origin and with c\`adl\`ag trajectories. See \cite{S99} for the basic theory of additive processes. 

Set ${\R}_0:={\R}-\{0\}.$ For any fixed $\epsilon>0,$ denote $S_{\epsilon}:=\{|x|>\epsilon\}\subseteq {\R}_0.$ Let us denote ${\cal B}$ and ${\cal B}_0$ the $\sigma-$algebras of Borel sets of $\R$ and ${\R}_0$ respectively. The distribution of an additive process can be caracterized by the triplet $(\Gamma_t, \sigma^2_t, \nu_t), \, t\geq 0$, where $\{\Gamma_t, t\geq 0\}$ is a continuous function null at the origin, $\{\sigma^2_t, t\geq 0\}$ is a continuous and non-decreasing function null at the origin and $\{\nu_t,t\geq 0\}$ is a set of L\'evy measures on ${\R}$, that is, a set of positive measures such that for any $t\geq 0,$ $\nu_t(\{0\})=0$ and $\int_{\R}(1\wedge x^2)\nu_t(dx)<\infty.$ Moreover, for any set $B\in {\cal B}_0$ such that $B\subseteq S_{\epsilon}$ for a certain $\epsilon>0$, $\nu_{\cdot}(B)$ is a continuous and increasing function null at the origin. 

If in addition we assume stationarity of the increments (namely, $X$ is a L\'evy process), then, for any $t\geq 0,$ the triplet becomes $(\gamma_L t, \sigma^2_L t, \nu_L t)$, where $\gamma_L$ is a real constant, $\sigma^2_L$ is a positive constant and $\nu_L$ is a L\'evy measure on $\R.$ 
Note that, thanks to the stationarity of the increments, a L\'evy process is fully characterized just by the triplet $(\gamma_L, \sigma^2_L, \nu_L),$ that is, the triplet in the case $t=1.$

Set $\Theta:=[0,\infty)\times \R.$ Let us denote by $\theta:=(t,x)$ the elements of $\Theta.$ Accordingly, $d\theta$ will denote the pair $(dt,dx).$ For $T\geq 0,$ we can introduce the measurable spaces $(\Theta_{T,\epsilon}, {\cal B}(\Theta_{T,\epsilon}))$ where $\Theta_{T,\epsilon}:=[0,T]\times S_{\epsilon}$ and ${\cal B}(\Theta_{T,\epsilon})$ is the corresponding Borel $\sigma$-field. Observe that $\Theta_{\infty, 0}=[0,\infty)\times {\R}_0$ and that $\Theta$ can be represented as $\Theta=\Theta_{\infty,0}\cup([0,\infty)\times\{0\}).$ Also observe that $[0,\infty)\times \{0\}\simeq [0,\infty).$

We can introduce a measure $\nu$ on $\Theta_{\infty,0}$ such that for any $B\in {\cal B}_0$ we have $\nu([0,t]\times B):=\nu_t(B).$ The hypotheses on $\nu_t$ guarantee that $\nu(\{t\}\times B)=0$ for any $t\geq 0$ and for any $B\in{\cal B}_0.$ In particular $\nu$ is $\sigma-$finite. 
Given $G\in {\cal B}(\Theta_{\infty, 0})$ we introduce the jump measure $N$ associated to $X$, defined as 

$$N(G)=\#\{t:\, (t,\Delta X_t)\in G\},$$
with $\Delta X_t=X_t-X_{t-}.$ Recall that $N$ is a Poisson random measure on ${\cal B}(\Theta_{\infty,0})$ with

$${\E}[N(G)]=  {\E}\big[ (N(G) - \E [N(G)] \,)^2 \big] =\nu(G).$$
Let $\widetilde N(dt,dx):=N(dt,dx)-\nu(dt, dx)$ be the compensated measure.

According to the L\'{e}vy-It\^{o} decomposition (see \cite{S99}) we can write:

\begin{equation}
X_t=\Gamma_t+W_t+J_t, \quad t\geq 0. \label{rep}
\end{equation}
Here $\Gamma$ is a continuous deterministic function null at the origin and 
$W$ is a centered Gaussian process with variance process $\sigma^2$ independent of $J$ (and $N$).
In relation with $W$ we can also define a $\sigma-$finite measure $\sigma$ on $[0,\infty)$ such that $\sigma([0,t])=\sigma^2_t.$
The process $J$ is an additive process with triplet $(0,0,\nu_t)$ defined by

\begin{equation}\label{purejumpdec}
J_t=\int_{\Theta_{t,1}} x N(ds,dx)+\lim_{\epsilon \downarrow 0}\int_{\Theta_{t,\epsilon}-\Theta_{t,1}} x {\widetilde N}(ds,dx),
\end{equation}
where the convergence is $ a.s.$ and uniform with respect to $t$ on every bounded interval. Following the literature, we will call the process $J=\{J_t, t\geq 0\}$ a pure jump additive process.

Moreover, if $\{{\cal F}^W_t, \,t\geq 0\}$ and $\{{\cal F}^J_t, \, t\geq 0\}$ are, respectively, the completed natural filtrations of $W$ and $J$, then, for every $t\geq 0,$ we have ${\cal F}^X_t={\cal F}^W_t \vee {\cal F}^J_t.$ The proof is the same as in the L\'evy case (see \cite{SUV07}).

We can consider on $\Theta$ the $\sigma-$finite Borel measure 

$$\mu(dt,dx):=\sigma(dt)\delta_0(dx)+\nu(dt,dx).$$ 
So, for $E\in{\cal B}(\Theta)$,

$$\mu(E)=\int_{E(0)}\sigma(dt)+\iint_{E'} \nu(dt, dx),$$ 
where $E(0)=\{t\geq 0:\, (t,0)\in E\}$ and $E'=E-E(0).$ Note that $\mu$ is continuous in the sense that $\mu(\{t\}\times B)=0$ for all $t\geq 0$ and $B\in {\cal B}.$ See \cite{D07} for a discussion on the importance of this condition for random measures with infinitely divisible distribution.
Then, for $E\in {\cal B}(\Theta)$ with $\mu(E)<\infty,$ we can define the measure

$$M(E)=\int_{E(0)}dW_t+{L^2-}\lim_{n\uparrow \infty}
\iint_{\{(t,x)\in E : \frac{1}{n}<|x|<n\}} \widetilde N(dt,dx),$$
that is a centered random measure with independent values such that ${\E}\big[M(E_1)M(E_2)]=\mu(E_1\cap E_2)$ for $E_1,E_2\in {\cal B}(\Theta)$ with 
$\mu(E_1) < \infty$ and $\mu(E_2) < \infty.$ The measure $M$ appears as a mixture of independent Gaussian and compensated Poisson random measures. We can write

$$M(dt,dx)=(W\otimes \delta_0)(dt, dx)+{\tilde N}(dt, dx).$$

\begin{remark}\label{casos}
\mbox{}

\begin{enumerate}
\item
If we take $\sigma^2\equiv 0$, $\mu=\nu$ and $M={\tilde N},$ we recover the Poisson random measure case.

\item
If we take $\nu=0$, we have $\mu(dt,dx)=\sigma(dt)\delta_0(dx)$ and $M(dt,dx)=(W\otimes \delta_0)(dt, dx)$ and we recover the independent increment centered Gaussian measure case. 

\item
If we take $\sigma^2_t:=\sigma_L^2 t$ and $\nu(dt,dx)=dt\nu_L(dx)$, we obtain $M(dt,dx)=\sigma_L (W\otimes \delta_0)(dt, dx)+{\tilde N}(dt, dx)$ and we recover the L\'evy case (stationary increments case).

\item
If $\nu=0$ and $\sigma^2_t=\sigma^2_L t$, we have $\mu(dt, dx)=\sigma^2_L dt\delta_0 (dx)$ and $M(dt,dx)=\sigma_L W(dt)\delta_0 (dx)$ and we recover the Brownian motion case.

\item
If $\sigma^2\equiv 0$ and $\nu(dt,dx)=dt\delta_1(dx)$, we have $\mu(dt,dx)= dt\delta_1(dx)$ and $M(dt,dx)={\tilde N}(dt)\delta_1(dx)$ and we recover the standard Poisson case.
\end{enumerate}
\end{remark}

\begin{remark}\label{tildemeasure}
A similar situation can be developed with the mixtures ${\bar \mu}(dt, dx)=\sigma(dt)\delta_0 (dx)+x^2\nu(dt, dx)$ and ${\bar M}=(W\otimes \delta_0) + x{\tilde N}$ as can be seen in \cite{I56} and \cite{SUV07} in the particular context of L\'evy processes.
\end{remark}

\begin{remark}
Given $\mu$ we can consider the Hilbert space $H:=L^2(\Theta, {\cal B}, \mu)$ and introduce the so-called isonormal additive process on $(\Omega, {\cal F}^X, {\prob}),$ i.e. a process $L:=\{L(h), h\in H\}$ such that $L$ is linear and 

$${\E}(e^{izL(h)})=\exp(\phi(z,h)), \, z\in {\R},$$
with 

$$\phi(z,h)=\int_{\Theta}((e^{izh(t,x)}-1-izh(t,x)){1\!\!1}_{{\R}_0}-\frac{z^2}{2}h^2(t,x){1\!\!1}_{\{0\}})\mu(dt,dx).$$
Observe that we can rewrite 

$$\phi(z,h)=\int_{\Theta_{\infty,0}}((e^{izh(t,x)}-1-izh(t,x))\nu(dt,dx)-\int_0^{\infty}\frac{z^2}{2}h^2(t,0)\sigma(dt).$$
See \cite{Y08} for the details. In the case $\nu\equiv 0$, $L$ becomes an isonormal Gaussian process (see \cite{N06}). Note also that ${\cal F}^L$, which is the natural completed $\sigma-$algebra generated by $L$, coincides with ${\cal F}^X$ and $M(A)=L({1\!\!1}_A)$, for any $A\in {\cal B}$.
\end{remark}

\section{Malliavin-Skorohod calculus for additive processes in $L^2.$}

\hspace{0.5cm} Here we summarize the Malliavin-Skorohod calculus with respect to the random measure $M$ on its canonical space in the $L^2-$framework. The construction is the same as for the stationary case and follows \cite{Y08}, but in a way it is close to \cite{SUV07} and \cite{ALPV08}. This is the first step towards our final goal of extending the calculus to the $L^1$ and $L^0$ frameworks.

\subsection{The chaos representation property}
\hspace{0.5cm} Given $\mu$, we can consider the spaces ${\Le}^2_n:=L^2\Big(\Theta^n, {{\cal B}(\Theta)}^{\otimes n},\mu^{\otimes n}\Big)$ and define the It\^{o} multiple stochastic integrals $I_n(f)$ with respect to $M$ for functions $f$ in ${\Le}^2_n$ by linearity and continuity starting from $I_n(f):=M(E_1)\cdots M(E_n)$
if $f:={1\!\!1}_{E_1\times\cdots\times E_n}$ with $E_1,\dots, E_n \in {\cal B}(\Theta)$ pairwise disjoint and with finite measure $\mu.$  In particular, for any $f\in {\Le}^2_n$ we have $I_n(f)=I_n({\tilde f})$, where $\tilde f$ is the symmetrization of $f$.
By construction, $I_n$ does not charge the diagonal sets, i.e. the sets 

$$\{(\theta_1,\dots, \theta_n)\in \Theta^n: \theta_{i_1}=\cdots =\theta_{i_k} \mbox{ for some different } i_1,\dots, i_k\in \{1,\dots,n\}\}.$$
So, we can consider $f\in {\Le}^2_{n}$ to be null on the diagonal sets.  
Then we have the so-called chaos representation property, that is, for any functional $F\in L^2(\Omega, {\cal F}^X, {\prob}),$ we have

$$F=\sum_{n=0}^\infty I_n(f_n)$$
for a certain unique family of symmetric kernels $f_n \in {\Le}^2_n.$ See \cite{I56} for details of this construction. 
For the chaos representation property see also 
Theorem 3.3 in \cite{DR} and Theorem 2.2 in \cite{D08}.

\subsection{The Malliavin and Skorohod operators}

\hspace{0.5cm} The chaos representation property of $L^2(\Omega,{\cal F}^X, {\prob})$ shows that this space has a Fock space
structure. Thus it is possible to apply all the machinery related to the anhilation operator (Malliavin derivative) and the
creation operator (Skorohod integral) as it is exposed, for example, in \cite{NV90}.

Consider $F=\sum_{n=0}^\infty I_n(f_n),$ with $f_n$ symmetric and such that $\sum_{n=1}^\infty n\, n ! \Vert f_n\Vert^2_{{\Le}^2_n} <\infty.$  The Malliavin derivative of $F$ is an object of $L^2(\Theta \times \Omega,\mu \otimes \prob)$, defined as

\begin{equation}
D_{\theta} F:=\sum_{n=1}^\infty n I_{n-1} \Big(f_n\big(\theta,\cdot \big)\Big),\ \ \theta\in \Theta. \label{derivada}
\end{equation}
We will denote by  ${\rm Dom}D$ the domain of this operator.

On other hand, let $u\in L^2\big(\Theta\times \Omega, {\cal B}(\Theta) \otimes {\cal F}^X, \mu\otimes \prob)$. For every
$\theta \in \Theta$ we have the chaos decomposition

$$u_{\theta}=\sum_{n=0}^\infty I_n(f_n(\theta,\cdot))$$
where $f_n \in {\Le}^2_{n+1}$ is symmetric in the last $n$ variables. Let ${\tilde f}_n$ be the symmetrization in all $n+1$ variables. Then we define the Skorohod integral of $u$ by

\begin{equation}
\delta(u):=\sum_{n=0}^\infty I_{n+1}(\tilde f_n),
\label{integral}
\end{equation}
in $L^2(\Omega),$ provided $u\in \text{Dom}\,\delta,$ that means $\sum_{n=0}^\infty (n+1)!\,\Vert \tilde f_n\Vert^2_{{\Le}^2_{n+1}}<\infty.$
Moreover if $u\in {\rm Dom}\,\delta$ and $F\in {\rm Dom}\, D$ we have the duality relation

\begin{equation}
{\E}[\delta(u)\,F]={\E}\int_{\Theta} u_{\theta} \, D_{\theta} F\,
\mu(d\theta). \label{dualitat}
\end{equation}

We recall that if $u\in Dom\delta$ is actually predictable with respect to the filtration generated by $X$, then the Skorohod integral coincides with the 
(non anticipating) It\^o integral in the $L^2-$setting with respect to $M.$

\subsection{The Clark-Haussmann-Ocone formula}

\hspace{0.5cm} Given $A\in {\cal B}(\Theta)$ we can consider the $\sigma-$algebra ${\cal F}_A$ generated by $\{M(A'): A'\in {\cal B}(\Theta),
A'\subseteq A\}.$ Following \cite{NV90} we have that $F$ is ${\cal F}_A-$measurable, if for any $n\geq 1,$ $f_n (\theta_1,\dots,\theta_n)=0,$ $\mu^{\otimes n} - a.e.$ unless $\theta_i\in A\quad \forall \,\, i =1,\dots n.$

In particular, we are interested in the case $A:=[0,t)\times {\R}.$ Let us denote from now on, ${\cal F}_{t-}:={\cal F}_{[0,t)\times \R}.$ Obviously, if $F\in {\rm Dom}\, D$ and it is ${\cal F}_{{t-}}-$measurable then $D_{s,x} F=0$ for a.e. $s\geq t$ and any $x\in {\R}.$
From the chaos representation property we can see that for $F\in L^2(\Omega),$

$${\E}[F|{\cal F}_{{t-}}]=\sum_{n=0}^{\infty} I_n\Big(f_n (\theta_1,\dots, \theta_n)\prod_{i=1}^n {1\!\!1}_{[0,t)} (t_i)\Big),$$ (see e.g. \cite{DOP08}). Then, for $F\in Dom D$, we have

$$D_{s,x} {\E}[F|{\cal F}_{t-}]= {\E}[D_{s,x}F|{\cal F}_{t-}] {1\!\!1}_{[0,t)}(s), \, (s,x)\in \Theta.$$

Using these facts and following Theorems 4.1, 12.16 and 12.20 of \cite{DOP08} (or the same steps as in Proposition 1.3.14 in \cite{N06}), we can prove the so-called Clark-Hausmann-Ocone (CHO) formula:

\begin{theorem}\label{CHO2 formula}

If $F\in Dom D$ we have $$F={\E}(F)+\delta({\E}[D_{t,x}F|{\cal F}_{t-}]).$$
\end{theorem}

\begin{remark}
Being the integrand a predictable process, the Skorohod integral $\delta$ in Theorem \ref{CHO2 formula} above is actually an It\^o integral. 
\end{remark}

\begin{remark}\label{decompactified2}
The CHO formula can be rewritten in a decompactified form as 

$$F={\E}(F)+\int_0^{\infty} {\E}(D_{s,0}F|{\cal F}_{s-})dW_s+\int_{\Theta_{\infty,0}} {\E}(D_{s,x}F|{\cal F}_{s-}){\tilde N}(ds,dx).$$
See \cite{BDLOP}.

\end{remark}

\section{A canonical space for additive processes}

\hspace{0.5cm} First we consider the pure jump case (process $J$) and then the general case (process $X$).

\subsection{A canonical space for $J$}

\hspace{0.5cm} We will set our work on the canonical space for $J$, introduced in \cite{SUV07}. 
Hereafter, we review the construction in a slightly different way, more convenient for our purposes, and in the more general context of additive processes. First we will consider the process on $\Theta_{T,\epsilon}$, for fixed $T>0$ and $\epsilon>0$, and then we will consider $\Theta_{\infty,0}$ taking $T\uparrow \infty$ and $\epsilon\downarrow 0.$

Assume for the moment that $\nu$ is concentrated on $\Theta_{T,\epsilon}$ or otherwise let us consider $\nu=\nu {1\!\!1}_{\Theta_{T, \epsilon}}.$ Observe that in particular, $\nu(\Theta_{T, \epsilon})<\infty.$ 
Note that in this case, 

$$c_{\epsilon}(t):=\int_0 ^t \int_{\epsilon<|x|\leq 1} x\nu(ds, dx)=\int_0^t \int_{-1}^1 x \nu(ds, dx)$$
and $|c_{\epsilon}(t)|\leq \nu(\Theta_{t, \epsilon})$ for any $t\in [0,T].$ Then, taking the characterization (\ref{purejumpdec}) into account, the process $J_t+c_{\epsilon}(t)$ can be identified with a time inhomogeneous compound Poisson process with parameter $\nu(\Theta_{T, \epsilon})$, that in particular has a finite number of jumps on $[0,T].$

Any trajectory of $J$ can be described by a finite sequence $\big((t_1,x_1),\dots, (t_n,x_n)\big)$, for some $n$, where $t_1,\dots, t_n\in [0,T]$: $t_1<t_2<\cdots <t_n$, are the jump instants and $x_1,\dots,x_n\in S_{\epsilon}$ are the corresponding sizes. Let $\alpha$ denote the empty sequence. 
 So we can define

\begin{enumerate}[(i)]

\item
$$\Omega^J_{T,\epsilon}:=\bigcup_{n\ge 0}\Theta_{T,\epsilon}^n,$$
where $\Theta_{T,\epsilon}^0=\{\alpha\}.$ Note that for any $n\neq n'$, $\Theta_{T,\epsilon}^n\cap \Theta_{T,\epsilon}^{n'}=\emptyset$;

\item
${\cal F}_{T,\epsilon}:=\sigma\big\{B\subset
\Omega^J_{T,\epsilon}:\, B=\bigcup_{n\ge 0} B_n\, (\mbox{disjoint}), \, B_n\in {\cal
B}\big(\Theta_{T,\epsilon}\big)^{\otimes n}\big\}=\bigvee_{n\ge 0}
{\cal B}(\Theta_{T,\epsilon})^{\otimes n}$;

\item
The probability measure 
$\prob_{T,\epsilon}$ such that, for $B=\bigcup_n B_n$ (pairwise disjoint) with $B_n\in {\cal
B}(\Theta_{T,\epsilon})^{\otimes n},$

$$\prob_{T,\epsilon} (B):=e^{-\nu (\Theta_{T,\epsilon})}\, \sum_{n=0}^\infty
\dfrac{\nu^{\otimes n}\big(B_n\big)}{n!},$$ where
$\nu^0=\delta_{\alpha}$.
\end{enumerate}

The pure jump process $\{J_t,\, t\in[0,T]\}$ on $(\Omega^J_{T,\epsilon}, {\cal F}_{T,\epsilon}, {\prob}_{T,\epsilon})$ is given by

\begin{equation}
J_t(\omega)=
\begin{cases}
\sum_{j=1}^nx_j\, {1\!\!1}_{[0,t]}(t_j)-\int_0^t\int_{-1}^1 x \nu(ds, dx),&  \text{if}\
\omega=\big((t_1,x_1),
\dots,(t_n,x_n)\big),\\
-\int_0^t\int_{-1}^1 x \nu(ds, dx), & \text{if}\   \omega=\alpha.
\end{cases} \label{canonic}
\end{equation}

Recall that given a measurable space $(E,{\cal E}),$ it is easy to see that the family of sets ${\cal E}^{\otimes n}_\text{sym}=\{C\in {\cal E}^{\otimes n}: \, C \ \text{is symmetric}\}$ is a $\sigma$-field. Here $C$ is symmetric if for all permutations $\pi$ of $\{1,\dots, n\}$ we have $C=\pi(C)=\{\pi(x): x\in C\}$ where $\pi(x):=(x_{\pi(1)}, \dots, x_{\pi(n)}).$ Recall also that a function $f:E^n\longrightarrow\R$ is ${\cal E}^{\otimes n}_\text{sym}$-measurable if and only if $f$ is ${\cal E}^{\otimes n}$-measurable and symmetric.
Let now ${\cal F}_{T,\epsilon, \text{sym}}$ be the sub-$\sigma$-field of ${\cal F}_{T,\epsilon}$ defined as 

$${\cal F}_{T,\, \epsilon,\,\text{sym}}:=\bigvee_{n\ge 0} {\cal B}\big(\Theta_{T,\epsilon}\big)^{\otimes n}_\text{sym}.$$ 
Let ${\cal F}^J_{T,\epsilon}$ be the $\sigma$-field generated by $J.$ It is easy to see that ${\cal F}^J_{T,\epsilon}={\cal F}_{T,\epsilon, \, {\rm sym}}.$
\vspace{0.5cm}

Now we extend the construction given above to the space $\Theta_{\infty,0}$ through a projective system of probability spaces.

First of all observe that $\Omega^J_{T,\epsilon}$ is a metric space. In fact for $u,v\in \Omega^J_{T,\epsilon}$, $u\in \Theta_{T,\epsilon}^n$, $v\in
\Theta_{T,\epsilon}^m$, we can define the distance

$$d(u,v):=\begin{cases}
1, & \text{ if }\ n\ne m, \ \text{ or }\ n=m\ \text{ and }\ d_{2n}(u,v)>1,\\
d_{2n}(u,v), & \text{ if }\  n=m\ \text{ and }\ d_{2n}(u,v)\le 1,
\end{cases}$$
where $d_k$ is the Euclidean distance on $\R^k$. Then $\Omega^J_{T,\epsilon}$ is a Polish space (metric, separable and
complete) and the $\sigma$-field ${\cal F}_{T,\epsilon}$ coincides with the Borel $\sigma$-field. We say that $(\Omega^J_{T,\epsilon},
{\cal F}_{T,\epsilon})$ is a separable standard Borel space. See Definition 2.2 in \cite{P65}.

For $m\ge 1$, let $(\Omega^J_m,{\cal F}_m,{\prob}_m):=(\Omega^J_{m,\frac{1}{m}},{\cal F}_{m,\frac{1}{m}},{\prob}_{m,\frac{1}{m}})$ be the canonical space corresponding to $\Theta_m:=[0,m]\times S_{\frac{1}{m}}.$
Observe that:

\begin{enumerate}

\item
$\{\Theta_m, m\geq 1\}$ and $\{\Omega^J_m, m\geq 1\}$ are increasing sequences of sets;

\item
$\Theta_{\infty, 0}=\cup_{m\geq 1} \Theta_m$ is an increasing union of sets;

\item
$\Theta_{\infty, 0}=\cup_{m\geq 1} (\Theta_m-\Theta_{m-1})$ is the union of pairwise disjoint sets. Remark that for $m=0$ we have the empty set.
\end{enumerate}

Consider the maps $\pi_m\colon \Omega^J_{m+1}\longrightarrow \Omega^J_m$ defined by

$$\pi_m\big((t_1,x_1),\dots,(t_r,x_r)\big)=\big((t_{i_1},x_{i_1}),\dots,
(t_{i_s},x_{i_s})\big),$$ where $(t_{i_1},x_{i_1})\dots,(t_{i_s},x_{i_s})$ are the points of
$(t_1,x_1)\dots, (t_r,x_r)$ belonging to $\Theta_m.$ If there are no points on this subspace we have
$\pi_m\big((t_1,x_1),\dots,(t_r,x_r)\big)=\alpha.$ It is straightforward to check that

$$\prob_m(B)=\prob_{m+1}(\pi_m^{-1}(B)),\quad \forall B\in{\cal F}_m.$$

The canonical space $\Omega^J$ for the pure jump additive process $J$ on $\Theta_{\infty, 0}$ can be defined as the projective limit of
the system $(\Omega^J_m,\pi_m, m\ge 1).$ Let ${\cal F}$ be the $\sigma$-field generated by the canonical projections $\overline
\pi_m:\Omega^J\to\Omega^J_m.$ Then, from \cite{P65}, there is a unique probability $\prob$ on $(\Omega^J,{\cal F})$, such that

$$\prob_m(B)=\prob(\overline\pi_m^{-1}(B)),\quad \forall B\in{\cal F}_m.$$
By construction, the projective limit $\Omega^J$ is the set of all sequences $(\omega^{(1)},\omega^{(2)},\dots,)$ with $\omega^{(m)}\in \Omega^J_m$
such that $\pi_{m}(\omega^{(m+1)})=\omega^{(m)}.$ In our setup, $\Omega^J=\cup_{n=0}^{\infty} \Theta^n_{\infty,0}$ and the probability measure $\prob$ is concentrated on the subset of $\Omega^J$ given by the following elements:   

\begin{itemize}
\item The empty sequence $\alpha$, corresponding to the element $(\alpha,\alpha,\dots)$.

\item All infinite sequences of pairs $(t_i,x_i)$ that are constant in the tail, that is, it exists $r>0$ such that $(t_{r+i},x_{r+i})=(t_r,x_r)$ for any $i\geq 0.$ This corresponds to the elements $(\omega^{(1)},\omega^{(2)},\dots,)$ such that $\omega^{(r)}=\omega^{(r+1)},\dots$ for some $r$. In this case we will usually write only the relevant finite part.

\item All  infinite sequences $((t_1,x_1),(t_2,x_2),\dots \big)$ such that for every $m>0$ there is only a finite number of $(t_i, x_i)$ on $\Theta_m.$
\end{itemize}
Furthermore, given the interpretation of $\Omega^J$ as the set of finite or infinite sequences 

$$\big((t_1,x_1),(t_2,x_2),\dots\big)$$ 
above exposed, the canonical projection
\begin{align*}
\overline \pi_m&\colon \Omega^J\longrightarrow\Omega^J_m\\
(\omega^{(1)},& \,\omega^{(2)},\cdots)\to \omega^{(m)}
\end{align*}
gives $\overline \pi_m\big((t_1,x_1),(t_2,x_2),\dots\big)$, which is then the finite sequence of points $(t_i,x_i)$ such that $t_i\in[0,m]$ and $|x_i|>\frac{1}{m}$. In the sequel, both $\Omega^J$ and $\overline\pi_m$ should be understood in this sense.

Now define the $\sigma$-field on $\Omega^J:$

$${\cal F}_{\text{sym}}:=\bigvee_{n\ge 0}\overline \pi_m^{-1}\big({\cal F}_{m,\, \text{sym}}\big).$$

Finally, the process $\{J_t,\, t\ge 0\}$ on $(\Omega^J,{\cal F}_{\text{sym}}, \prob)$ can be defined as follows. 
For any $t$, if
$\omega=(\omega^{(m)})_{m\ge 1}\in \Omega^J,$ set

$$J_t(\omega)=\lim_n\sum_{m=1}^n\big( J^{(m)}_t(\omega^{(m)})-J^{(m-1)}_t(\omega^{(m-1)})\big)$$
assuming $J^{(0)}_t\equiv 0.$ Here the convergence is ${\prob}-$a.s. and the $J^{(m)}$ are given as in
(\ref{canonic}). Moreover ${\cal F}^J$, the $\sigma-$algebra generated by $J,$ is equal to ${\cal F}_{\text{sym}}$.
The existence of the limit above is proved exactly as the It\^{o}-L\'evy representation of a pure jump L\'evy process which
gives the convergence a.s., uniform on $t\in[0,T]$, for any $T>0$, of an equivalent sequence.
Moreover, computing
the characteristic function, it is straightforward to see that $J=\{J_t,\, t\ge 0\}$ is a c\`adl\`ag additive process with triplet $(0,0,\nu_t).$
Observe that, in general, a random variable $F$ on $\Omega^J$ is given as

$$F(\omega)=\mbox{a.s.}-\lim_n\sum_{m=1}^n \big( F({\bar \pi}_m (\omega))-F({\bar \pi}_{m-1}(\omega))\big)+F(\alpha)=\mbox{a.s.}-\lim_m
F({\bar \pi}_m (\omega)),$$ provided these limits exists.

\subsection{A canonical space for $X$}

\hspace{0.5cm} Let $(\Omega^W,{\cal F}^W, {\prob}^W_{\sigma})$ be the canonical Wiener space and $\{\overline W^{\sigma}_t,\, t\geq 0\})$ be the canonical centered Gaussian process with independent increments and variance process $\sigma^2.$ That is, $\Omega^W=C_0 ([0,\infty)),$ is the space of continuous functions on $[0,\infty)$, null at the origin, with the topology of the uniform convergence on the compacts, ${\cal F}^W$ is the Borel $\sigma-$algebra and ${\prob}^W_{\sigma}$ is the probability measure that makes the projections ${\bar W}^{\sigma}_t:\Omega^W\longrightarrow \R$, $t\geq 0,$ be a centered Gaussian process (with independent increments) with variance process $\sigma^2.$
Let $(\Omega^J,{\cal F}^J, {\prob}^J, \{\overline J_t,t\geq 0\})$ be the canonical pure jump additive process associated to the measure $\nu$ defined before. 
We consider the product  space $(\Omega^W\times \Omega^J, {\cal F}^W\otimes {\cal F}^J,{\prob}^W_{\sigma}\otimes {\prob}^J)$ and put $W_t(\omega,\omega'):=\overline W^{\sigma}_t(\omega)$ and $J_t(\omega,\omega'):=\overline J_t(\omega').$ 
Finally, we consider the continuous deterministic function $\Gamma.$ 
Then 

$$X_t=\Gamma_t+W_t+J_t$$
is the canonical additive process with triplet $(\Gamma_t, \sigma^2_t, \nu_t), \, t\geq 0.$

\section{A Malliavin-Skorohod type calculus for $J$ on the canonical space}

\hspace{0.5cm} In this section we establish the operators and the basic calculus rules of a Malliavin-Skorohod calculus with respect to a pure jump additive process on the canonical space.

\subsection{An abstract duality relation}

\hspace{0.5cm} Let $\theta=(s,x)\in \Theta_{\infty,0}.$ Let $\omega\in\Omega^J$, that is, 
$\omega:=(\theta_1,\dots,\theta_n,\dots),$ with $\theta_i:=(s_i,x_i).$
We introduce the following two transformations from $\Theta_{\infty,0}\times \Omega^J$ to $\Omega^J:$

$$\epsilon^+_{\theta}\omega:=\big((s,x),(s_1,x_1),(s_2,x_2),\dots \big) ,$$
where a jump of size $x$ is added at time $s$, and 

$$\epsilon^-_{\theta}\omega:=\big((s_1,x_1),(s_2,x_2),\dots \big)-\{(s,x)\} ,$$
where we take away the point $\theta = (s,x)$ from $\omega$.

Observe that $\epsilon^+$ is well defined on $\Omega^J$ except on the set $\{(\theta,\omega): \theta\in \omega\}$, which has null $\nu\otimes {\prob}$ measure. 
We can set by convention that on this set, $\epsilon^+_{\theta}\omega:=\omega.$ The case of $\epsilon^-_{\theta}$ is also clear. In fact this operator satisfies $\epsilon^-_{\theta}\omega=\omega$ except on the set $\{(\theta,\omega): \theta\in \omega\}.$ 
For simplicity of the notation, when needed, we will denote ${\hat\omega}_i:=\epsilon^-_{\theta_i}\omega.$

These two transformations are analogous to the ones introduced in \cite{P96a}, where they are called creation and annihilation operators. Some of the results presented here have their correspondent in that paper, but our proofs are constructive on the canonical space. This differs from the approach used in \cite{P96a} and extends substantially the ideas presented in \cite{NV95}, from the standard Poisson to the additive case. See \cite{M93} for general information about creation and annihilation operators in quantum probability. 

Let $L^0(\Omega^J)$ denote the set of random variables defined on $\Omega^J$ and by $L^0 (\Theta_{\infty,0}\times \Omega^J)$ the set of measurable stochastic processes 
defined on $\Theta_{\infty,0}\times \Omega^J.$ Now we consider the following two definitions:

\begin{definition}
For a random variable $F\in L^0(\Omega^J)$, we define the operator 

$$T: L^0 (\Omega^J) \longmapsto L^0 (\Theta_{\infty,0}\times \Omega^J),$$
such that $(T_{\theta}F)(\omega):=F(\epsilon_{\theta}^+ \omega).$
\end{definition}

If $F$ is a ${\cal F}^J$-measurable, then

$$(T_{\cdot} F)(\cdot)\colon \Theta_{\infty, 0}\times \Omega^J \longrightarrow
\R$$ is ${\cal B}(\Theta_{\infty,0}) \otimes {\cal F}^J-$
measurable. Moreover, it $F=0,\,{\prob}\text{-a.s.}$, then $T_{\cdot} F(\cdot)=0,\, \nu\otimes {\prob}\text{-a.e.}$
So, $T$ is a closed linear operator defined on the entire $L^0 (\Omega^J).$ 
See \cite{SUV07} for a proof.

\vspace{2mm}
If we want to secure $T_{\cdot}F(\cdot)\in L^1(\Theta_{\infty,0}\times \Omega^J)$, we have to restrict the domain and guarantee that 

$$\E\int_{\Theta_{\infty,0}} |T_{\theta} F| \nu(d\theta)<\infty.$$
Remark that this requires a condition that is strictly stronger than $F\in L^1 (\Omega^J).$ Concretely, we have to assume that 

$$\sum_{m=1}^{\infty} e^{-\nu(\Theta_m-\Theta_{m-1})}\sum_{n=0}^{\infty} \frac{n}{n!} \int_{(\Theta_m-\Theta_{m-1})^n} |F(\theta_1,\dots,\theta_n)|\nu(d\theta_1)\dots\nu(d\theta_n)<\infty,$$
whereas $F\in L^1(\Omega)$ is equivalent only to 

$$\sum_{m=1}^{\infty} e^{-\nu(\Theta_m-\Theta_{m-1})}\sum_{n=0}^{\infty} \frac{1}{n!} \int_{(\Theta_m-\Theta_{m-1})^n} |F(\theta_1,\dots,\theta_n)|\nu(d\theta_1)\dots\nu(d\theta_n)<\infty.$$

\begin{definition}
For a random field $u\in L^0(\Theta_{\infty, 0}\times \Omega^J)$ we define the operator 

$$S: Dom S\subseteq L^0(\Theta_{\infty,0}\times\Omega^J)\longrightarrow L^0 (\Omega^J)$$
such that 

$$(Su)(\omega):=\int_{\Theta_{\infty,0}} u_{\theta}({\epsilon^-_{\theta}\omega})N(d\theta,\omega):=\sum_{i} u_{\theta_i}({\hat \omega}_i)<\infty.$$
In particular, if $\omega=\alpha$, we define $(Su)(\alpha):=0.$
\end{definition}

The operator $S$ is well defined on $L^1 (\Theta_{\infty,0}\times \Omega^J)$ as the following 
proposition says: 

\begin{proposition}\label{elau}
If $u\in L^1(\Theta_{\infty,0}\times \Omega^J),$ $Su$ is well defined and takes values in $L^1(\Omega).$ Moreover

$$\E\int_{\Theta_{\infty,0}} u_{\theta}({\epsilon^-_{\theta}\omega})N(d\theta,\omega)=\E\int_{\Theta_{\infty,0}} u_{\theta}(\omega)\nu(d\theta).$$
\end{proposition}

\noindent
\begin{proof}
Fix $\Omega^J_m$ and denote, for any $n\geq 0,$ $\omega:=(\theta_1,\dots,\theta_n)$ and $\theta:=(s,x).$ Denote also $c_m:=e^{-\nu(\Theta_m)}.$
We have

\begin{eqnarray*}
&&\E({1\!\!1}_{\Omega^J_m}\int_{\Theta_m} u_{\theta}({\epsilon^-_{\theta}\omega})N(d\theta,\omega))\\
&=&\sum_{n=1}^{\infty}\frac{c_m}{n!} \int_{\Theta^n_m}
\sum_{i=1}^n u_{\theta_i}(\theta_1,\dots,{\hat \theta}_i,\dots,\theta_n)\nu(d\theta_1)\cdots \nu(d\theta_n)\\
&=&\sum_{n=1}^{\infty}\frac{c_m}{n!} \int_{\Theta^n_m}
n u_{\theta}(\theta_1,\dots, \theta_{n-1})\nu(d\theta_1)\cdots \nu(d\theta_{n-1}) \nu(d\theta)\\
&=&\sum_{n=1}^{\infty}\frac{c_m}{(n-1)!} \int_{\Theta^{n-1}_m}\int_{\Theta_m}
u_{\theta}(\theta_1,\dots, \theta_{n-1}) \nu(d\theta_1)\cdots \nu(d\theta_{n-1}) \nu(d\theta)\\
&=&\sum_{l=0}^{\infty}\frac{c_m}{l!} \int_{\Theta^{l}_m}\int_{\Theta_m}
u_{\theta}(\theta_1,\dots, \theta_l) \nu(d\theta_1)\cdots \nu(d\theta_l) \nu(d\theta)\\
&=&\E({1\!\!1}_{\Omega^J_m}\int_{\Theta_m} u_{\theta}\nu(d\theta))
\end{eqnarray*}
The general case comes from dominated convergence. 
\end{proof}
\vspace{0.5cm}

\begin{remark}
We have proved that $L^1(\Theta_{\infty,0}\times\Omega^J)\subseteq Dom S.$ Moreover $S$ is closed in $L^1$ as an operator from $L^1(\Theta_{\infty,0}\times\Omega^J)$ to $L^1(\Omega).$ In fact, if we take a sequence $u^{(n)}\in L^1(\Theta_{\infty,0}\times\Omega^J)$ converging to $0$ in this space and we assume that $Su^{(n)}$ converges to $G$ in $L^1(\Omega^J)$, we can show that $G=0$. This is immediate because 

$$\E|G|\leq \E|G-Su^{(n)}|+\E|Su^{(n)}|.$$ 
Moreover, the first term in the right hand side converges to $0$ by hypothesis and the second one, using Proposition \ref{elau}, can be bounded by 
 
$$\E|Su^{(n)}|\leq \E\int_{\Theta_{\infty,0}} |u^{(n)}_{\theta} (\epsilon^-_{\theta}\omega)|N(d\theta,\omega)
=\E\int_{\Theta_{\infty,0}} |u^{(n)}_{\theta}|\nu(d\theta),$$
which also converges to $0$ by hypothesis. 
\end{remark}

\begin{remark}\label{predict}
Given $\theta=(s,x)$, for any $\omega$, we can define ${\tilde \omega}_{s}$ as the restriction of $\omega$ to jump instants strictly before $s.$ In this case, obviously, $\epsilon^-_{\theta}{\tilde \omega}_s={\tilde \omega}_s.$ 
If $u$ is predictable we have $u_{\theta}(\omega)
=u_{\theta}({\tilde \omega}_s)$. In this case, we have

$$u_{\theta}({\epsilon^-_{\theta}\omega})=u_{\theta}(\tilde{(\epsilon^-_{\theta}\omega)_s})=
u_{\theta}({\tilde\omega}_s)=u_{\theta}(\omega),$$ 
and   

$$(Su)(\omega)=\int_{\Theta_{\infty,0}} u_{\theta}({\epsilon^-_{\theta}\omega})N(d\theta,\omega)=\int_{\Theta_{\infty,0}} u_{\theta}(\omega)N(d\theta,\omega).$$
\end{remark}

Hereafter we introduce a fundamental relationship between the two operators $S$ and $T$:

\begin{theorem}\label{primeradualitat}
Consider $F\in L^0 (\Omega^J)$ and $u\in Dom S.$ Then $F\cdot Su\in L^1 (\Omega^J)$ if and only if 
$TF \cdot u\in L^1(\Theta_{\infty, 0}\times \Omega^J)$ and in this case

$$\E(F\cdot Su)=\E\int_{\Theta_{\infty,0}} T_{\theta} F \cdot u_{\theta}\,
\nu(d\theta).$$
\end{theorem}

\noindent
\begin{proof}
Using the fact that $F$ is symmetric, i.e. ${\cal F}_{sym}-$measurable in the canonical space, we have

\begin{eqnarray*}
&&\E(F\cdot Su \cdot {1\!\!1}_{\Omega^J_m})\\
&=&\sum_{n=1}^{\infty}\frac{c_m}{n!}\int_{{\Theta}_m^n}
F(\theta_1,\dots \theta_n) (Su)(\theta_1,\dots,
\theta_n) \nu(d\theta_1)\cdots \nu(d\theta_n)\\
&=&\sum_{n=1}^{\infty}\frac{c_m}{n!}\int_{{\Theta}_m^n} F(\theta_1,\dots
\theta_n) \sum_{i=1}^n u_{\theta_i}({\hat\omega}_i)
 \nu(d\theta_1)\cdots \nu(d\theta_n)\\
&=&\sum_{n=1}^{\infty}\sum_{i=1}^n \frac{c_m}{n!}\int_{{\Theta}_m^n}
T_{\theta_i}F(\theta_1,\dots, {\hat \theta}_i, \dots \theta_n)
u_{\theta_i}({\hat \omega}_i) \nu(d\theta_1)\cdots \nu(d\theta_n)\\
&=&\sum_{n=1}^{\infty} n \frac{c_m}{n!}\int_{{\Theta}_m^{n-1}}\int_{\Theta_m}
T_{\theta} F(\theta_1,\dots \theta_{n-1}) u_{\theta}({\hat \omega}_n)
 \nu(d\theta_1)\cdots \nu (d\theta_{n-1})\nu(d\theta)\\
&=&\E\Big({1\!\!1}_{\Omega^J_m}\int_{\Theta_m} T_{\theta}F u_\theta \,\nu(d\theta)\Big)
\end{eqnarray*}
Finally, we extend the result to $\Omega^J$ using the dominated
convergence theorem.
\end{proof}

\bigskip\noindent
Moreover we obtain the following rules of calculus:

\begin{proposition}\label{nova1}
If $u$ and $TF\cdot u$ belong to $DomS$, then we have $F\cdot Su=S(TF\cdot u),$   $\prob-a.e.$
\end{proposition}

\noindent
\begin{proof}
This is an immediate consequence of the fact that $T_{\theta_i}F({\hat\omega}_i)=F(\omega).$
\end{proof}

\begin{proposition}\label{nova2}
If $u$ and $Tu$ are in $Dom S$, then $T_{\theta}(Su)=u_{\theta}+S(T_{\theta}u),$ $\nu\otimes\prob-a.e.$
\end{proposition}

\noindent
\begin{proof}
For the left-hand side term we have

$$T_{\theta}(Su)(\omega) =(Su)(\epsilon^+_{\theta}\omega)=u_{\theta}(\omega)+\sum_{i} u_{\theta_i}(\epsilon^-_{\theta_i}\epsilon^+_{\theta}\omega)$$
and for the right-hand side term we have 

$$u_{\theta}(\omega)+S(T_{\theta}u)(\omega) =u_{\theta}(\omega) +\sum_{i} u_{\theta_i}(\epsilon^+_{\theta}\epsilon^-_{\theta_i}\omega).$$
The equality comes from 
$\epsilon^-_{\theta_i}\epsilon^+_{\theta}\omega=\epsilon^+_{\theta}\epsilon^-_{\theta_i}\omega,$ $\nu\otimes\prob-a.e.$
\end{proof}

\subsection{The intrinsic gradient and divergence operators and their duality}

\hspace{0.5cm} 
With the results of the previous sections we are ready to introduce two operators which also turn out to fullfill a duality relationship. These operators will be hereafter called intrinsic operators, being defined constructively on the canonical space.

We define the operator $$\Psi_{\theta}:=T_{\theta}-Id.$$ 
Observe that this operator is linear,
closed and satisfies the property

$$\Psi_{\theta}(F\, G)= G \,\Psi_{\theta} F+ F\, \Psi_{\theta}G+\Psi_{\theta}(F)\,
\Psi_{\theta}(G).$$

On other hand, for $u\in L^0 (\Theta_{\infty,0}\times \Omega^J)$ we consider the operator:

$${\cal E}: Dom {\cal E}\subseteq L^0 (\Theta_{\infty,0}\times \Omega^J)\longrightarrow L^0 (\Omega^J)$$
such that 

$$({\cal E}u)(\omega):=\int_{\Theta_{\infty, 0}} u_{\theta}(\omega)\nu(d\theta).$$
Note that $Dom{\cal E}$ is the subset of processes in $L^0 (\Theta_{\infty,0}\times \Omega^J)$ such that $u(\cdot,\omega)\in L^1(\Theta_{\infty,0})$, ${\prob}-$a.s. On other hand recall that, for $\omega$ fixed, we have $\epsilon^-_{\theta}\omega=\omega$, if $\theta\ne \theta_i$ for any $i$, and that $\nu(\{\theta: \theta=\theta_i, \mbox{ for some } i\})=0.$ So,  

\begin{equation}\label{igualtatnu}
\int_{\Theta_{\infty, 0}} u_{\theta}({\epsilon^-_{\theta}\omega})\nu(d\theta)=\int_{\Theta_{\infty, 0}} u_{\theta}(\omega)\nu(d\theta), \quad {\prob}-a.s.
\end{equation}
Then, for $u\in {Dom\Phi}:=Dom S\cap Dom {\cal E}\subseteq L^0(\Theta_{\infty,0}\times \Omega^J)$, we define

$$\Phi u:=Su-{\cal E}u.$$

\begin{remark}
Observe that $L^1(\Theta_{\infty,0}\times\Omega^J)\subseteq Dom\Phi.$
\end{remark}

\begin{remark}
Observe that from Proposition \ref{elau} and (\ref{igualtatnu}) we have that $E(\Phi u)=0,$ for any $u\in L^1 (\Theta_{\infty,0}\times \Omega).$
\end{remark}

\begin{remark}\label{predict2}
From Remark \ref{predict} and (\ref{igualtatnu}) we have

$$\Phi(u)=\int_{\Theta_{\infty,0}} u_{\theta}(\omega){\tilde N}(d\theta,\omega),$$
for any predictable $u\in Dom \Phi$. 
\end{remark}
As a corollary of Theorem \ref{primeradualitat} we have the following result:

\begin{proposition}\label{dual1}
Consider $F\in L^0 (\Omega^J)$ and $u\in Dom\Phi.$ Assume also $F\cdot u\in L^1(\Theta_{\infty,0}\times \Omega^J).$ Then $F\cdot {\Phi}u\in L^1 (\Omega^J)$ if and only if ${\Psi} F \cdot u\in L^1(\Theta_{\infty,0}\times \Omega^J)$ and in this case

$$\E(F\cdot {\Phi} u)=\E\Big(\int_{\Theta_{\infty,0}} {\Psi}_{\theta} F \cdot u_{\theta}\,
\nu(d\theta)\Big).$$
\end{proposition}

\vspace{2mm}
Analogously to the previous subsection we have also the following two results that can be proved immediately using Propositions \ref{nova1} and \ref{nova2} and recalling the definitions $\Psi=T-Id$ and $\Phi=S-{\cal E}.$

\begin{proposition}
\label{calc1}
If $F\in L^0(\Omega^J)$ and $u$, $F\cdot u$ and $\Psi F\cdot u$ belong to $Dom\Phi$, then we have 

$$F\cdot \Phi u=\Phi(F\cdot u)+\Phi(\Psi F\cdot u)+{\cal E}(\Psi F \cdot u), \quad\prob-a.s.$$
\end{proposition}

\begin{proposition}
\label{calc2}
If $u$ and $\Psi u$ belong to $Dom \Phi$, we have 

$$\Psi_{\theta}(\Phi u)=u_{\theta}+\Phi(\Psi_{\theta}u), \quad \nu\otimes\prob-a.e.$$
\end{proposition}

\begin{remark}\label{dotze}
If we change $\nu(ds,dx)$ by $x^2\nu(ds,dx)$ and we define the operators

$${\bar \Psi}_{s,x} F:=\frac{T_{s,x}F-F}{x},$$

$${\bar S} u (\omega):=\int_{\Theta_{\infty,0}}u_{s,x}(\epsilon^-_{s,x}\omega)x N(ds,dx),$$

$$({\bar {\cal E}}u)(\omega):=\int_{\Theta_{\infty, 0}} u_{s,x}(\omega)x^2 \nu(ds,dx)$$
and

$$\bar \Phi:={\bar S}-{\bar {\cal E}},$$
we can prove similar results to the previous ones. For example, if $F\in L^0 (\Omega^J)$, $u\in {Dom {\bar \Phi}},$ and $F\cdot u\in L^1(\Theta_{\infty,0}\times \Omega^J),$ then $F\cdot {\bar \Phi}u\in L^1 (\Omega^J)$ if and only if ${\bar \Psi} F \cdot u\in L^1(\Theta_{\infty,0}\times \Omega^J)$ and in this case

$$\E(F\cdot {\bar \Phi} u)=\E\Big(\int_{\Theta_{\infty,0}} {\bar \Psi}_{s,x} F
\cdot u_{s,x}\, x^2 \nu(ds, dx)\Big).$$

Note that the domains of $\bar \Psi$ and $\Psi$ are slightly different in view of the different measure $\nu.$ This has natural consequences also on the evaluations in $L^1.$ For example, 

$${\E}\int_{\Theta_{\infty,0}}|{\bar \Psi}_{s,x}F|x^2 \nu(ds, dx)={\E}\int_{\Theta_{\infty,0}} |{\Psi}_{s,x}F||x| \nu(ds, dx)\neq {\E}\int_{\Theta_{\infty,0}} |{\Psi}_{s,x}F| \nu(ds, dx).$$
\end{remark}

\subsection{Relationships between the intrinsic operators and the Malliavin-Skorohod operators.}

\hspace{0.5cm} 
In the last part of this section we study the intrinsic operators $\Psi$ and $\Phi$ in comparison with the Malliavin derivative and Skorohod integral defined in Section 3.2 restricted to the pure jump case, i.e. associated with ${\tilde N}(ds,dx)$. We will write $D^J$ and $\delta^J$, respectively. 
We show that the intrinsic operators are extensions of the two classical concepts.
 
First we need to recall some preliminary results.
The following key lemma is proved in \cite{SUV07} (see the proof of Lemma 5.2) and it is an
extension of Lemma 2 in \cite{NV95}.

\begin{lemma}
For any $n\geq 1$, consider the set

$$\Theta_{T,\epsilon}^{n,*}=\{(\theta_1,\dots, \theta_n)\in \Theta_{T,\epsilon}^{n}: \theta_i\neq \theta_j \mbox{ if } i\neq j\}.$$
Then, for any $g_k\in L^2 (\Theta^{k,*}_{\infty, 0})$, for $k\geq 1$, and $\omega\in \Omega^J$ we have

$$I_k (g_k)(\omega)=\int_{\Theta_{T,\epsilon}^{k,*}} g_k (\theta_1\dots,\theta_k) {\tilde N}(\omega, d\theta_1)\cdots {\tilde N}(\omega,d\theta_k), \quad \prob-a.e.$$
\end{lemma}

\noindent
\begin{proof}
Both expressions coincide for simple functions and define bounded linear operators. We remark that $g_k$ does not need to be symmetric.
\end{proof}

\vspace{2mm}
\noindent
The relationships between $D^J$ and $\Psi$, and $\delta^J$ and $\Phi$ are given by the following results, which extend corresponding results for the standard Poisson process given in \cite{NV95}.

\begin{lemma}\label{derivadasimple}
For a fixed $k\geq 0,$ consider $F=I_k (g_k)$ with $g_k$ a symmetric function of $L^2(\Theta_{\infty,0}^{k,*})$. Then, $F$ belongs to $Dom D^J\cap Dom \Psi$ and 

$$D^J I_k (g_k)=\Psi I_k (g_k), \quad \nu \otimes {\prob}-\text{a.e.}$$
\end{lemma}

\noindent
\begin{proof}
The fact that $F\in Dom D^J\cap Dom \Psi$ is obvious. From the definition of $\Psi$ we obtain

$$\Psi_{\theta} I_k (g_k)(\omega)= I_k (g_k)(\epsilon^+_{\theta}\omega)-I_k (g_k)(\omega)$$

$$=\int_{\Theta_{\infty,0}^{k,*}} g_k (\theta_1,\dots, \theta_k) {\tilde N}(\epsilon^+_{\theta}\omega, d\theta_1)
\cdots {\tilde N}(\epsilon^+_{\theta}\omega, d\theta_k)-\int_{\Theta_{\infty,0}^{k,*}}
g_k (\theta_1,\dots, \theta_k) {\tilde N}(\omega, d\theta_1)
\cdots {\tilde N}(\omega, d\theta_k)$$

$$=\int_{\Theta_{\infty,0}^{k,*}} g_k (\theta_1,\dots, \theta_k) \prod_{i=1}^k
({\tilde N}(\omega, d\theta_i)+N(\theta, d\theta_i))-\int_{\Theta_{\infty,0}^{k,*}} g_k (\theta_1,\dots,
\theta_k) \prod_{i=1}^k {\tilde N}(\omega, d\theta_i).$$
Using the fact that $g_k$ is null on the diagonals, only the integrals with $k-1$ integrators of type $\tilde N$ and one integrator of type $N$ remain. 
Using the fact that $g_k$ is symmetric in the last expression we obtain

$$
\Psi_{\theta} I_k (g_k)(\omega)=
k\int_{\Theta_{\infty,0}^{k-1,*}} g_k (\theta_1,\dots, \theta_{k-1},\theta) {\tilde N}(\omega, d\theta_1)
\cdots {\tilde N}(\omega, d\theta_{k-1})=D^J_{\theta} I_k (g_k).$$
\end{proof}

\begin{lemma}\label{lemauk}
For fixed $k\geq 1,$ consider $u_{\theta}=I_k (g_k (\cdot, \theta))$ where $g_k (\cdot, \cdot)\in L^2 (\Theta_{\infty,0}^{k+1,*})$ is symmetric with
respect to the first $k$ variables. Assume also $u\in Dom\Phi.$ Then,

$$\Phi(u)=\delta^J (u), \quad {\prob}-\text{a.e.}.$$
\end{lemma}

\noindent
\begin{proof}
First of all, note that

$$\delta^J(I_k (g_k (\cdot, \theta))(\omega)=I_{k+1} ({\tilde g}_k (\cdot,\cdot))(\omega)=I_{k+1} (g_k (\cdot,\cdot))(\omega)$$

$$=\int_{\Theta_{\infty,0}^{k+1,*}} g_k (\theta_1,\dots,\theta_k,\theta){\tilde N}(\omega,d\theta_1)\cdots {\tilde N}(\omega,
d\theta_k)N(\omega, d\theta)-\int_{{\Theta}_{\infty, 0}} u_{\theta}\nu(d\theta)$$

$$=\sum_j \int_{\Theta_{\infty,0}^{k,*}} g_k (\theta_1,\dots, \theta_k, \theta^0_j){\tilde N}(\omega,
d\theta_1)\cdots {\tilde N}(\omega,d\theta_k)-\int_{{\Theta}_{\infty, 0}} u_{\theta} \nu(d\theta),$$ where
the different $\theta_j^0$ are the jump points of $\omega=(\theta^0_1,\theta^0_2,\dots).$

Recall that ${\tilde g}_k$, the symmetrization of $g_k$ with respect to all its variables, is null on the diagonals, so
$\theta^0_j$ has to be different of all $\theta_i$, for $i=1,\dots,k.$
Now observe that we can write ${\tilde N}(\omega,d\theta)=N(\theta_j^0, d\theta)+{\tilde N}({\epsilon^-_{{\theta}_j^0}\omega}, d\theta)$, where for simplicity we write ${\hat \omega}_j:=\epsilon^-_{\theta_j^0}\omega.$ Then we have, 

\begin{eqnarray*}
\delta^J (u)&=&\sum_{j} \sum_{l=0}^k \binom{k}{l}\int_{\Theta_{\infty,0}^{k,*}} g_k(\theta_1,\dots, \theta_k, \theta^0_j)N(\theta^0_j, d\theta_1)\cdots N(\theta^0_j, d\theta_l){\tilde N}({\hat \omega}_j, d\theta_{l+1})\cdots {\tilde N}({\hat \omega}_j,
d\theta_k)\\
&-&\int_{{\Theta}_{\infty, 0}} u_{\theta}\nu(d\theta)\\
&=&\sum_{j} \int_{\Theta_{\infty,0}^{k,*}} g_k (\theta_1,\dots, \theta_k, \theta^0_j){\tilde N}({\hat \omega}_j,
d\theta_{1})\cdots {\tilde N}({\hat \omega}_j,d\theta_k)-\int_{{\Theta}_{\infty, 0}} u_{\theta}\nu(d\theta)=\Phi(u).
\end{eqnarray*}
\end{proof}

\begin{remark}
Recall that $u\in L^2(\Theta_{\infty,0}\times\Omega^J)$ does not imply that $u\in L^1(\Theta_{\infty,0}\times\Omega^J)$, nor that $u\in Dom\Phi.$
\end{remark}

\begin{theorem}\label{derivada2}
Let $F\in L^2(\Omega^J)$. Then $F\in {Dom D}^J $ if and only if $ \Psi F\in L^2(\Theta_{\infty,0}\times \Omega^J)$ and in this case we have

$$D^J F=\Psi F, \quad \nu\otimes \prob-\text{a.e.}$$
\end{theorem}

\noindent
\begin{proof}
Note that $F\in Dom\Psi$ because $Dom\Psi$ is the entire $L^0(\Omega^J).$ Consider $u_{\theta}=I_k (g_k (\cdot,\theta))$ as in Lemma \ref{lemauk}, that is, we are assuming also that $u\in Dom \Phi$. Then from \eqref{dualitat}, Lemma \ref{lemauk} and Proposition \ref{dual1} we have formally that

\begin{equation}\label{id1}
\E\int_{\Theta_{\infty, 0}} D^J_{\theta} F u_{\theta}\nu(d\theta)=\E(F\delta^J (u))=\E(F\Phi (u))
=\E\int_{\Theta_{\infty,0}} \Psi_{\theta} F u_{\theta}\nu(d\theta).
\end{equation}
The objects in \eqref{id1} are well defined either if $F\in Dom D^J$ or if $\Psi F\in L^2(\Theta_{\infty,0}\times\Omega^J)$. 
In particular the previous equalities are true in the case $g_k(\theta_1,\dots, \theta_k,\theta):={1\!\!1}_{A_1}(\theta_1)\cdots{1\!\!1}_{A_k}(\theta_k){1\!\!1}_{A}(\theta)$
for any collection of pairwise disjoint and measurable sets $A_1,\dots, A_k, A$ with finite measure $\nu$. In fact in this case $u\in L^1(\Theta_{\infty,0}\times\Omega^J)\subseteq {Dom\Phi}.$
So, in particular we have 

$${\E}(I_k({1\!\!1}^{\otimes k}_{A_1\times\cdots\times A_k})\int_A D^J_{\theta}F\nu(d\theta))={\E}(I_k({1\!\!1}^{\otimes k}_{A_1\times\cdots\times A_k})\int_A \Psi_{\theta}F\nu(d\theta)).$$
By linearity and continuity we conclude that
$\Psi F=D^J F,$ $\nu\otimes {\prob}-a.e.$
\end{proof}

\begin{theorem}\label{integral2}

Let $u\in L^2(\Theta_{\infty, 0}\times \Omega^J)\cap Dom\Phi.$ Then $u\in {Dom \delta}^J$ if and only if $ \Phi u\in
L^2(\Omega^J)$ and in this case we have

$$\delta^J u=\Phi u, \quad {\prob}-\text{a.s.}$$
\end{theorem}

\noindent
\begin{proof} Let $G=I_k (g_k)$ as in Lemma \ref{derivadasimple}. Note that $G$ is in $Dom D^J.$ Then from \eqref{dualitat}, Lemma \ref{derivadasimple} and Proposition \ref{dual1} we have formally that

\begin{equation}\label{id2}
\E(\delta^J (u)G)=E\int_{\Theta_{\infty, 0}} u_{\theta} D^J_{\theta} G
\nu(d\theta)
=\E\int_{\Theta_{\infty, 0}} u_{\theta}
\Psi_{\theta} G \nu(d\theta)=E(G\Phi(u)).
\end{equation}
The  objects in \eqref{id2} are well defined if either $\Phi(u)\in L^2 (\Omega^J)$ or if $u\in Dom\delta^J$  hold. Then the conclusion follows.
\end{proof}

\begin{remark}
Similar results can be obtained for the operators $\bar \Phi$ and $\bar \Psi.$ See Remarks \ref{tildemeasure} and \ref{dotze}.
\end{remark}

\section{The Clark-Hausmann-Ocone formula}

\subsection{The CHO formula in the pure jump case}

\hspace{0.5cm} As an application of the previous results in the pure jump case we present a CHO type formula as an integral representation of random variables in $L^1(\Omega^J)$.
This in particular extends the formula proved in \cite{P96a} for the standard Poisson case, as well as the formulae of CHO type proved in the $L^2$ setting, 
see e.g. \cite{BDLOP}, \cite{DOP08}.

\begin{theorem}\label{jumpsCHO}
Let $F\in L^1 (\Omega^J)$ and assume  $\Psi F\in L^1(\Theta_{\infty,0}\times \Omega^J).$ Then 

$$F=\E(F)+\Phi (\E({\Psi}_{t,x} F|{\cal F}_{t-})) \qquad \prob- a.s.$$
\end{theorem}

\noindent
\begin{proof}
The argument is organised in two steps. 
\begin{enumerate}
\item
Assume first that we are in $\Omega^J_m.$ In this case, $\nu$ is a finite measure concentrated on $\Theta_m.$

Given $F\in L^1(\Omega^J)$ we can define, for every $n\geq 1,$ $F_n$ such that $F_n=F$ if $|F|\leq n$, $F_n=n$ if $F_n\geq n$ and $F_n=-n$ if $F\leq -n.$ Of course, $F_n\in L^2(\Omega^J).$ And moreover $|F_n|\leq |F|$ for any $n$ and 

$$|\Psi F_n|\leq |TF_n|+|F_n|\leq |TF|+|F|\leq |\Psi F|+2|F|.$$

Applying Theorems \ref{CHO2 formula}, \ref{derivada2} and \ref{integral2} we obtain

$$F_n={\E}(F_n)+\Phi(\E(\Psi_{\theta} F_n|{\cal F}_{t-})), \quad \prob- a.s.$$

Being $\nu$ finite, we note that $\Psi F_n\in L^2(\Theta_{\infty,0}\times \Omega^J)$ and $E(\Psi_{\theta} F_n|{\cal F}_{t-})\in L^2(\Theta_{\infty,0}\times \Omega^J)\cap Dom\Phi$.  

Using Remark \ref{predict2} we obtain

$$F_n={\E}(F_n)+\int_{\Theta_{m}}\E(\Psi_{\theta} F_n|{\cal F}_{t-}){\tilde N}(d\theta), \quad \prob- a.s.$$

Clearly, $F_n-{\E}(F_n)$ converges in $L^1$ to $F-{\E}(F).$ So, to prove the formula for $F\in L^1(\Omega^J_m)$ it is enough to prove that  

$$\int_{\Theta_{m}}\E(\Psi_{\theta} (F-F_n)|{\cal F}_{t-}){\tilde N}(d\theta)\longrightarrow_{n\uparrow \infty} 0,$$ with convergence in $L^1 (\Omega^J).$
Indeed we have 

$$|\int_{\Theta_{m}}\E(\Psi_{\theta} (F-F_n)|{\cal F}_{t-}){\tilde N}(d\theta)|\leq \int_{\Theta_{m}}\E(|\Psi_{\theta} (F-F_n)||{\cal F}_{t-})N(d\theta)+\int_{\Theta_{m}}\E(|\Psi_{\theta} (F-F_n)||{\cal F}_{t-})\nu(d\theta).$$

So, it is enough to show that both summands on the right-hand side converge to $0.$ Observe that using Proposition \ref{elau} and Remark \ref{predict} these 
two quantities have the same expectation, which is equal to 

$${\E}\int_{\Theta_{m}}|\Psi_{\theta} (F-F_n)|\nu(d\theta).$$

Now, the sequence $|\Psi (F-F_n)|$ converges to $0,$ $\prob$-a.s., and it is dominated by 

$$|\Psi (F-F_n)|\leq |\Psi F|+|\Psi F_n|\leq 2(|\Psi F|+|F|),$$
as this last quantity belongs to $L^1(\Theta_m\times\Omega^J_m)$ by hypothesis. 

\item
Now we consider the general case. Then we have 

$$F{1\!\!1}_{\Omega^J_m} - {\E}(F{1\!\!1}_{\Omega^J_m})={1\!\!1}_{\Omega^J_m}\int_{\Theta_m}\E(\Psi_{\theta} F|{\cal F}_{t-}){\tilde N}(d\theta), \quad \prob-a.s.$$

It is immediate to see that if $F\in L^1 (\Omega^J)$ the left-hand side of the equality converges to $F-{\E}(F).$ The convergence of the right-hand side is a consequence of the fact that 

$$  \E \Big( {1\!\!1}_{\Omega^J_m} \Big\vert \int_{\Theta_m} \E(\Psi_{\theta} F|{\cal F}_{t-} ) \tilde N(d\theta)\Big\vert \Big)
\leq 2 \int_{\Theta_{\infty,0}}{\E}(|\Psi_{\theta} F|)\nu(d\theta)$$
and the dominated convergence.      
\end{enumerate}
By this we end the proof.
\end{proof}

\begin{remark}
Observe that under the conditions of the previous theorem we have 

$$\Psi_{s,x}\E[F|{\cal F}_{{t-}}]=\E[\Psi_{s,x} F|{\cal F}_{{t-}}]{1\!\!1}_{[0,t)}, \quad \nu\otimes {\prob}-a.e.$$
Indeed, on $\Omega^J_m$ we consider the functionals $F_n$ introduced in the proof of the previous theorem and  we have 

$$\Psi_{s,x}\E[F_n|{\cal F}_{{t-}}]=\E[\Psi_{s,x} F_n|{\cal F}_{{t-}}]{1\!\!1}_{[0,t)}, \quad \nu\otimes {\prob}-a.e.$$
The sequence $\Psi_{s,x}F_n$ converges a.s. to $\Psi_{s,x}F$ and the term is bounded in $L^1(\Theta_{\infty,0}\times \Omega^J_m)$, so the right-hand side term converges to $\E[\Psi_{s,x} F|{\cal F}_{{t-}}]{1\!\!1}_{[0,t)}{1\!\!1}_{\Omega^J_m}.$ Then, the left-hand side has a limit in $L^1$. On other hand, this left-hand side term also converges $\nu\otimes {\prob}$-a.e. to $\Psi_{s,x}\E[F|{\cal F}_{{t-}}].$ So, the result follows.  
\end{remark}

\begin{example}
 
Consider a pure jump additive process $L$, i.e. for all $t$, $L_t$ can be represented by the following L\'evy-It\^o decomposition:

$$L_t=\Gamma_t+\int_0^t \int_{\{|x|>1\}} xN(ds,dx)+\int_0^t \int_{\{|x|\leq 1\}}x{\tilde N}(ds,dx).$$
Consider $L_T$ (for $T>0$). If we assume ${\E}(|L_T|)<\infty,$ or equivalently that 

$$
\int_0^T \int_{\{|x|>1\}} |x| \nu(ds,dx)<\infty
$$
(see \cite{CT} Proposition 3.13), then we can write 

$$L_T =\Gamma_T+\int_0^T \int_{\{|x|>1\}} x\nu(ds,dx)+\int_0^T \int_{\R}x{\tilde N}(ds,dx).$$

On the other hand, applying the CHO formula, we have 

$$\Psi_{s,x} L_T=x{1\!\!1}_{[0,T]}(s)$$ 
and 

$${\E}(\Psi_{s,x}L_T|{\cal F}_{s-})=x{1\!\!1}_{[0,T)}(s).$$ 
So, the conditions of Theorem 6.1 are equivalent to  

\begin{equation}
\label{finiteexp}
{\E}\int_0^T\int_{{\R} }|x| \nu(ds,dx)<\infty.
\end{equation}
So, under this condition, the CHO formula gives

$$L_T={\E}(L_T)+\int_0^T \int_{\R}x{\tilde N}(ds,dx).$$
This is clearly coherent with the L\'evy-It\^o decomposition because, under \eqref{finiteexp}, we have

$${\E}(L_T)=\Gamma_T+\int_0^T \int_{\{|x|>1\}} x\nu(ds,dx).$$
\end{example}

\begin{example}

Let $X:=\{X_t,t\in [0,T]\}$ be a pure jump L\'evy process with triplet $(\gamma_L t,0,\nu_L t).$ Let $S_t:=e^{X_t}$ be an asset price process (see e.g. \cite{CT} for the use of exponential L\'evy models in finance). Let $\Q$ be a risk-neutral measure. Recall that $e^{-rt}e^{X_t}$ is a $\Q-$martingale under the following assumptions on $\nu_L$ and $\gamma_L$: 

$$\int_{|x|\geq 1} e^x \nu_L(dx)<\infty$$ 
and 

$$\gamma_L=\int_{\R}(e^y-1-y{1\!\!1}_{\{|y|<1\}})\nu(dy).$$
See \cite{CT} or \cite{JV} for details.
These conditions allow us to write without lost of generality, 

$$X_t=x+(r-c_2)t+\int_0^t\int_{\R} y {\tilde N}(ds,dy),$$ 
where 

$$c_2:=\int_{\R} (e^y-1-y)\nu_L(dy)$$
and $N$ is a Poisson random measure under $\Q.$
According to Theorem 6.1, if $F=S_T\in L^1(\Omega^J)$ and ${\E}_{\Q}[\Psi_{s,x}S_T|{\cal F}_{s-}]\in L^1([0,T]\times \Omega^J)$ we have 

$$S_T={\E}_{\Q}(S_T)+\int_{\Theta_{T,0}} {\E}_{\Q}[\Psi_{s,x}S_T|{\cal F}_{s-}]{\tilde N}(ds,dx).$$ 
Observe that $\Psi_{s,x}S_T(\omega)=S_T (e^x-1), \quad \ell\times\nu_L\times {\Q}- a.s.,$ and this process belongs to $L^1( \Theta_{\infty,0}\times \Omega^J)$ if and only if $\int_{\R}|e^x-1|\nu_L(dx)<\infty.$
Here $\ell$ denotes the Lebesgue measure on $[0,T]$.
Then, in this case, we have 

$$S_T={\E}_{\Q}(S_T)+\int_{\Theta_{T,0}} e^{r(T-s)} (e^x-1)S_{s-}{\tilde N}(ds,dx).$$ 
So, this result covers L\'evy processes with finite activity and L\'evy processes with infinite activity but finite variation. 
\end{example}

\subsection{The CHO formula in the general case}

\hspace{0.5cm} For the sake of completeness we present a version of the CHO formula in the general additive case that extends the formula in 
Remark \ref{decompactified2} from the $L^2$ setting to the $L^1$ setting. 

Let $W$ be an isonormal Gaussian process indexed by a Hilbert space $H.$ 
The classical construction of the Malliavin derivative for functionals of an isonormal Gaussian process is as follows, see e.g. \cite{N06}. 
Let $\cal S$ be the space of smooth functionals of type $F=f(W(h_1),\dots, W(h_n))$ where $f\in C_b^{\infty}({\R}^n)$ and $h_1,\dots, h_n$ are elements of $H.$ For a given $F\in {\cal S},$ its Malliavin derivative is the $H-$valued random variable defined as  

$${\cal D}^W F:=\sum_{i=1}^n (\partial_i f)(W(h_1),\dots,W(h_n))h_i.$$
Associated to these definition and for any $p\geq 1$, we can define the space ${\D}^{1,p}$ as the closure of $\cal S$ with respect the norm

$$||F||_{1,p}:=(\E(|F|^p+||{\cal D}^W F||^p_H)^{\frac{1}{p}}.$$ 
In particular we can consider the spaces ${\D}^{1,2}$ and ${\D}^{1,1}$, as the closures with respect the norms  

$$||F||_{1,2}:=(\E(|F|^2+||{\cal D}^W F||^2_H))^{\frac{1}{2}}.$$ 
and 

$$||F||_{1,1}:=E(|F|)+E(||{\cal D}^W F||_H),$$ respectively.
Observe that we have the inclusions ${\D}^{1,p}\subseteq L^p (\Omega^W)$ and that ${\D}^{1,2}\subseteq {\D}^{1,1}.$
By closure, tha Malliavin derivative can be defined in any space ${\D}^{1,p}.$ 

In particular, if $H:=L^2([0,\infty),\sigma)$, the Gaussian process $W$ introduced in Section 2 is an isonormal Gaussian process on $H.$ Then, for any $F\in {\D}^{1,1}$, we have the following version of the CHO formula (see \cite{KOL}):

\begin{theorem}\label{CHOW}
For any $T>0$ and $F\in {\D}^{1,1}$ we have 

$$F=\E(F)+\int_0^T {\E}({\cal D}^W_tF|{\cal F}_{t-})dW_t \qquad \prob- a.s.$$
\end{theorem}

In this case we can also relate the operator ${\cal D}^W$ with the operator $D_{t,0}$, which is restricted to the Gaussian case (compare with \eqref{derivada}). 
%We now denote this operator $D^W$, similarly to the notation used for the pure jump case. 
We have also the following results (see \cite{N06}):

\begin{proposition}
Let $F\in L^2(\Omega^W)$ such that $F\in Dom {\cal D}^W.$ Then ${\cal D}^W F\in L^2( [0,\infty)\times \Omega^W)$ if and only if $ F\in Dom D_{t,0}$
and in this case, 

\begin{equation}
D_{t,0} F={\cal D}_{t}^{W}F.
\end{equation}
\end{proposition}

Recall now that $\Theta=\Theta_{\infty,0}\cup ([0,\infty)\times\{0\})$ and $\Omega=\Omega^W\times \Omega^J$, hence $\omega=(\omega^W,\omega^J)\in \Omega^W\times \Omega^J.$ 
Using the independence between $W$ and $J$ we interpret ${\cal D}^W$ and $\Psi$ as operators on $L^0(\Omega^W\times\Omega^J)\cong L^0(\Omega^W; L^0(\Omega^J))$ and $L^0(\Omega^W\times\Omega^J)\cong L^0(\Omega^J; L^0(\Omega^W))$ respectively, on their suitable domains. See \cite{SUV07} for a similar construction in $L^2(\Omega^W\times\Omega^J).$
Compare also with \cite{DOP08}.

Now, for $F\in L^0 (\Omega^W\times \Omega^J)$ we define the operator

\begin{equation}\label{deffor}
{\nabla}_{t,x}F:={1\!\!1}_{\{0\}}(x){\cal D}^W_t F+{1\!\!1}_{{\R}_0}(x) \Psi_{t,x} F
\end{equation} 
on the domain

$$Dom {\nabla}:={\D}^{1,1}(\Omega^W; L^0(\Omega^J))\cap L^0 (\Omega^J; L^0(\Omega^W)).$$
Note that $\nabla$ extends $D_{t,x}$ from ${\D}^{1,2}(\Omega)$ to $Dom\nabla.$ Note also that in the right-hand side of \eqref{deffor}, if $\sigma\equiv 0$ only the second term remains and if $\nu\equiv 0$ only the first term remains. 

Then, we have the following result

\begin{corollary}
If $F\in L^2 (\Omega)\cap Dom\nabla,$ we have 

$$\Psi F\in L^2 (\Theta_{\infty,0}\times \Omega) \mbox{ and } {\cal D}F\in L^2([0,\infty)\times \Omega)\Longleftrightarrow
F\in Dom D,$$ and in this case

\begin{equation}\label{derbarm}
D_{t,x}F={\nabla}_{t,x}F \qquad \mu\times {\prob}-a.e.
\end{equation}
\end{corollary}

Hence we can extend the CHO formula to the following theorem:

\begin{theorem}\label{globalCHO}
Let $F\in L^1 (\Omega)\cap Dom \nabla$ and assume $\Psi F\in L^1(\Theta_{\infty,0}\times\Omega).$
Then, 

$$F=\E(F)+\int_{\Theta_{\infty,0}} {\E}({\Psi}_{s,x} F|{\cal F}_{s-}){\tilde N}(ds,dx)+ 
\int_0^{\infty} {\E}({\cal D}^W_s F|{\cal F}_{s-})dW_s^{\sigma}\qquad {\prob}-a.s.$$
\end{theorem}

\noindent
\begin{proof}
The result can be proved applying Remark \ref{decompactified2} to the approximating sequence $F_n$ introduced in Theorem \ref{jumpsCHO} and using Theorems \ref{jumpsCHO} and \ref{CHOW}.
\end{proof}

\begin{remark}
This CHO formula identifies the kernels of the predictable representation property proved in Theorem 8 in \cite{C13}, in the case of additive integrators. 
\end{remark}

\section{Integration with respect pure jump volatility modulated Volterra processes}

\hspace{0.5cm} Consider a pure jump volatility modulated additive driven Volterra ($\cal{VMAV}$) process $X.$ 
The definition of a $\cal {VMAV}$ process is the extension of the definition of a pure jump volatility modulated L\'evy driven Volterra ($\cal{VMLV}$) process as described in \cite{BBPV}. The process $X$ is given as 

\begin{equation}\label{VMLV}
X(t)=\int_0^t g(t,s)\sigma(s)dJ(s)
\end{equation}   
provided the integral is well defined. Here $J$ is a pure jump additive processes, $g$ is a deterministic function and $\sigma$ is a predictable process with respect the natural completed filtration of $J.$ 

Recall that using the L\'evy-It\^o representation $J$ can be written as 

$$J(t)=\Gamma_t+\int_{\Theta_{t,0}-\Theta_{t,1}} x{\tilde N}(ds,dx)+\int_{\Theta_{t,1}} x N(ds,dx),$$
where $\Gamma$ is a continuous deterministic function that we assume of bounded variation in order to admit integration with respect $d\Gamma.$ 
Recall also that in the case $\int_0^t \int_{|x|>1} |x|\nu(ds,dx)<\infty,$ we can rewrite the previous expression as 
$$J(t)=\Gamma_t+\int_{\Theta_{t,0}} x{\tilde N}(ds,dx)+\int_{\Theta_{t,1}} x \nu(ds, dx).$$
For each $t$, the integral \eqref{VMLV} is well defined (see \cite{BBPV}) if the following hypotheses are satisfied: 

\begin{eqnarray*}
(H1)&\quad& \int_{0}^{\infty} |g(t,s)\sigma(s)|d\Gamma_s<\infty,\\
(H2)&\quad& \int_{\Theta_{\infty,0}} 1\wedge \big(g(t,s)\sigma(s)x\big)^2\nu(dx,ds)<\infty,\\
(H3)&\quad& \int_{\Theta_{\infty,0}} \big|g(t,s)\sigma(s)x \big[{1\!\!1}_{\{|g(t,s)\sigma(s)x|\leq 1\}}-{1\!\!1}_{\{|x|\leq 1\}}\big] \big|\nu(dx,ds) <\infty.
\end{eqnarray*}

Hereafter we discuss the problem of defining an integral with respect to $X$ as integrator, i.e. to give a meaning to
$$
\int_0^t Y(s) dX(s)
$$
for a fixed $t$ and a suitable stochastic processes $Y$.

Indeed, exploiting the representation of $J$, an integration with respect to $X$ can be treated as the sum of integrals with respect to the corresponding components of 
$J.$ That is, it is enough to define integrals with respect $\int_0^t g(t,s)\sigma(s)d\Gamma_s$, $\int_0^t \int_{|x|\leq 1} g(t,s)\sigma(s)x{\tilde N}(ds,dx)$ and 
$\int_0^t \int_{|x|>1} g(t,s)\sigma(s)x N(ds,dx).$ Under the assumption that $\Gamma$ has finite variation and using the fact that $N$ is of finite variation on $\{|x|>\delta\}$, for any $\delta>0$, the integration with respect to the first and third term presents no difficulties. We have to discuss the second term, specifically the case when $J$ has infinite activity and infinite variation and the corresponding $X$ is not a semimartingale. In fact, if $X$ was a semimartingale, we could perform the integration in the It\^o sense. 
However, in general, $X$ is not a semimartingale. We can refer to \cite{BBPV} for the characterization of the restrictions on $g$ to guarantee the semimartingale structure of $X$. Also in \cite{BBPV} a definition of an integral with respect to a non semimartingale $X$ driven by a L\'evy process is given by means of  the Malliavin-Skorohod calculus. Their technique is naturally constrained to an $L^2$ setting.

Within the framework presented in this paper, we can extend the definition proposed in \cite{BBPV} to reach out for additive noises beyond the $L^2$ setting. Specifically we can present the following result:

\begin{theorem}\label{condicions}
Assume the following hypothesis on $X$ and $Y$:
\begin{enumerate}
\item
For $s\geq 0$, the mapping $t\longrightarrow g(t,s)$ is of bounded variation on any interval $[u,v]\subseteq (s,\infty).$
\item
The function 
$${\cal K}_g(Y)(t,s):=Y(s)g(t,s)+\int_s^t \big(Y(u)-Y(s) \big) \, g(du,s), \quad t > s,$$
is well defined a.s., in the sense that $Y(u)-Y(s)$ is integrable with respect to $g(du,s)$ as a pathwise Lebesgue-Stieltjes integral.  
\item
The mappings 
$$(s,x)\longrightarrow {\cal K}_g(Y)(t,s)\sigma(s)x{1\!\!1}_{\Theta_{t,0}-\Theta_{t,1}}(s,x)$$
and 
$$(s,x)\longrightarrow \Psi_{s,x}({\cal K}_g(Y)(t,s)\sigma(s))x{1\!\!1}_{\Theta_{t,0}-\Theta_{t,1}}(s,x)$$
belong to $Dom\Phi.$
\end{enumerate}
Then, the following integral, is well defined: 
\begin{eqnarray*}
\int_0^t Y(s) d(\int_0^s \int_{|x|\leq 1} g(s,u)\sigma(u)x{\tilde N}(du,dx))&:=&\Phi(x{\cal K}_g(Y)(t,s)\sigma(s){1\!\!1}_{\Theta_{t,0}-\Theta_{t,1}}(s,x))\\
&+&\Phi(x\Psi_{s,x}({\cal K}_g(Y)(t,s))\sigma(s){1\!\!1}_{\Theta_{t,0}-\Theta_{t,1}}(s,x))\\
&+&{\cal E}(x\Psi_{s,x}({\cal K}_g(Y)(t,s))\sigma(s){1\!\!1}_{\Theta_{t,0}-\Theta_{t,1}}(s,x)).
\end{eqnarray*}
\end{theorem}
The result is proved following the same lines given in \cite{BBPV}. 
The proof relies on the definitions of $\Phi$, $\Psi$ and the calculus rules of Propositions \ref{calc1} and \ref{calc2}.

Here we stress that the theory presented in Section 5 of this paper allows to go beyond Definition 3 in \cite{BBPV} and to treat pure jump additive processes $J$, in particular without second moment.
To be specific, in the finite activity case, $L^2(\Theta_{\infty,0}\times \Omega^J) \subseteq L^1(\Theta_{\infty,0}\times \Omega^J)$. 
Then Theorem \ref{condicions} is an extension of Definition 3 in \cite{BBPV}. In particular, for example, hypothesis (3) in Theorem \ref{condicions} is verified if the two mappings are in $L^1(\Theta_{\infty,0}\times\Omega^J)$ for any $t\geq 0.$ 
On the contrary, in the infinite activity case, Theorem \ref{condicions} reaches cases not covered by Definition 3 in \cite{BBPV} and viceversa. 

\bigskip
As illustration we give an example of a pure jump L\'evy process without second moment as driver $J$ and we consider a kernel function $g$ of shift type, i.e. it only depends on the difference $(t-s)$. The chosen kernel appears in applications to turbulence, see \cite{BBPV}, \cite{BS03} and references therein. 

\begin{example}
Assume $L$ to be a symmetric $\alpha-$stable L\'evy process, for $\alpha\in (0,2)$, see e.g. \cite{CT}, corresponding to the triplet $(0,0,\nu_L)$ with $\nu_L(dx)=c|x|^{-1-\alpha}dx.$ Recall that in the case $\alpha\leq 1$ the process is of finite variation whereas if $\alpha>1,$ the process is of infinite variation. Take 
$$g(t,s):=(t-s)^{\beta-1}e^{-\lambda(t-s)}{1\!\!1}_{[0,t)}(s)$$
with $\beta\in (0,1)$ and $\lambda>0.$ 
Note that 
$$g(du,s)=-g(u,s)(\frac{1-\beta}{u-s}+\lambda)du.$$  
Take $\sigma\equiv 1.$ 
We concentrate on the component 
$$J(t)=\int_{\Theta_{t,0}-\Theta_{t,1}} x{\tilde N}(ds,dx),$$
and so on the definition of the integral 
\begin{equation}\label{mainintegral}
X(t):=\int_0^t g(t,s)dJ(s)=\int_0^t \int_{|x|\leq 1} g(t,s)x {\tilde N}(ds,dx),
\end{equation}
for $t \geq 0$. 
As anticipated, the component
$$\int_{\Theta_{t,1}}xN(ds,dx)$$ 
is of finite variation and the corresponding integral 
$$\int_0^t\int_{|x|>1} g(t,s)x N(ds,dx)=\sum_{i} g(t,s_i)x_i$$ 
presents no problems because the sum is $\prob$-a.s. finite. 
In relation with \eqref{mainintegral}, which is not a semimartingale (see \cite{BP}), we have four situations:
\begin{enumerate}
\item
If $\alpha\in (0,1)$ and $\beta>\frac{1}{2}$, $g(t,s)x$ belongs to $L^1\cap L^2$ 
\item
If $\alpha\in [1,2)$ and $\beta>\frac{1}{2}$, $g(t,s)x$ belongs to $L^2$ but not to $L^1.$
\item
If $\alpha\in (0,1)$ and $\beta\leq\frac{1}{2}$, $g(t,s)x$ belongs to $L^1$ but not to $L^2.$
\item
If $\alpha\in [1,2)$ and $\beta\leq \frac{1}{2}$, $g(t,s)x$ belongs not to $L^2$ nor to $L^1.$
\end{enumerate}
The case (1) is both covered by Definition 3 in \cite{BBPV} and our Theorem \ref{condicions}. 
The case (2) is covered by \cite{BBPV}, while the case (3) is only covered by our Theorem \ref{condicions}.
It seems not possible to cover case (4). Recall that the domain of $\Psi$ is $L^0(\Omega)$, but the domain of $\Phi$ is only slightly greater than $L^1(\Theta_{\infty,0}\times \Omega).$

\vspace{2mm}
Just to show the types of computation involved, let us consider the particular case of a  $\mathcal{VMAV}$ process as integrand. Namely,
$$
Y(s)=\int_0^s \int_{{|x|\leq 1}} \phi(s-u)x {\tilde N}(du,dx), \quad 0 \leq s \leq t,
$$ 
where $\phi$ is a positive continuous function such that the integral $Y$ is well defined. Consider the case $\alpha<1$ and $\beta\in (0,1).$ In order to see that $\int_0^t Y(s-)dX(s)$ is well defined we have to check:
\begin{enumerate}
\item
The process $Y(u)-Y(s)$ is integrable with respect to $g(du,s)$ on $(s,t]$, as a pathwise Lebesgue-Stieltjes integral.
\item
The mappings 
$$
(s,x)\longrightarrow x{\cal K}_g(Y)(t,s){1\!\!1}_{[0,t]}(s){1\!\!1}_{\{|x|\leq 1\}}$$
and 
$$(s,x)\longrightarrow x\Psi_{s,x}({\cal K}_g(Y)(t,s)){1\!\!1}_{[0,t]}(s){1\!\!1}_{\{|x|\leq 1\}}
$$
belong to $Dom\Phi.$
\end{enumerate}
We have
\begin{eqnarray*}
{\cal K}_g(Y)(t,s)&=&g(t,s)\int_{[0,s)}\int_{|x|\leq 1} \phi(s-v)x {\tilde N}(dv,dx)\\
&-&\int_s^t g(u,s)(\frac{1-\beta}{u-s}+\lambda)\int_{[s,u)} \int_{|x|\leq 1}\phi(u-v)x{\tilde N}(dv,dx)du\\
&-&\int_s^t g(u,s)(\frac{1-\beta}{u-s}+\lambda)\int_{[0,s)}\int_{|x|\leq 1}[\phi(u-v)-\phi(s-v)]x{\tilde N}(dv,dx)du.
\end{eqnarray*}
In terms of $\Phi$ we can rewrite  
\begin{eqnarray*}
{\cal K}_g(X)(t,s)&=&g(t,s)\Phi(\phi(s-\cdot)x{1\!\!1}_{\{|x|\leq 1\}}{1\!\!1}_{[0,s)})\\
&-&\int_s^t g(u,s)(\frac{1-\beta}{u-s}+\lambda)\Phi(\phi(u-\cdot)x{1\!\!1}_{\{|x|\leq 1\}}{1\!\!1}_{[s,u)}(\cdot))du\\
&-&\int_s^t g(u,s)(\frac{1-\beta}{u-s}+\lambda)\Phi([\phi(u-\cdot)-\phi(s-\cdot)]x{1\!\!1}_{\{|x|\leq 1\}}{1\!\!1}_{[0,s)}(\cdot))du.
\end{eqnarray*}
Moreoever, using Proposition \ref{calc2}, we have 
$$
\Psi_{s,x}{\cal K}_g (X)(t,s)=-x{1\!\!1}_{\{|x|\leq 1\}}\int_s^t g(u,s)\phi(u-s)(\frac{1-\beta}{u-s}+\lambda){1\!\!1}_{[0,u)}(s)du.
$$
So, it is enough to check that the two mappings 
$$
(s,x)\longrightarrow x{\cal K}_g(Y)(t,s){1\!\!1}_{[0,t]}(s){1\!\!1}_{\{|x|\leq 1\}}
$$
and 
$$(s,x)\longrightarrow x\Psi_{s,x}({\cal K}_g(Y)(t,s)){1\!\!1}_{[0,t]}(s){1\!\!1}_{\{|x|\leq 1\}}$$
are in $L^1(\Theta_{\infty,0}\times \Omega).$ 

To proceed further with the illustration we consider the case $\phi(y)=y^{\gamma}$ with $\gamma>0$ and $\beta+\gamma\geq 1.$
Note that using Proposition \ref{elau} we have  
$${\E}|\int_a^b \int_{|x|\leq 1} f(v)x{\tilde N}(dv,dx)|\leq c\int_a^b |f(v)|dv,$$
where $0\leq a\leq b$, $c$ is a generic constant and $f$ is an integrable function. \\
We study the first mapping. 
For the first term of this first mapping we have  
$${\E}\int_0^t\int_{|y|\leq 1} |y|\, g(t,s)\, \Big|\int_{[0,s)}\int_{|x|\leq 1} \phi(s-v)x {\tilde N}(dv,dx)\Big| \nu(dy)ds$$
$$\leq c \int_0^t g(t,s)\int_0^s \phi(s-v)dvds\leq B(\beta,\gamma+1)t^{\beta+\gamma}<\infty.$$
Recall that $$B(m,n):=\int_0^1 t^{m-1} (1-t)^{n-1}dt \quad (m,n>0).$$
The second term satisfies
$${\E}\int_0^t\int_{|y|\leq 1}|y|\int_s^t g(u,s)(\frac{1-\beta}{u-s}+\lambda) \Big|\int_{[s,u)} \int_{|x|\leq 1}\phi(u-v)x{\tilde N}(dv,dx) \Big|du \nu(dy)ds$$
$$\leq c\int_0^t\int_s^t g(u,s)(\frac{1-\beta}{u-s}+\lambda)\int_s^u \phi(u-v)dvduds$$
$$\leq c\int_0^t\int_s^t (u-s)^{\beta+\gamma}(\frac{1-\beta}{u-s}+\lambda)duds<\infty.$$
Finally for the third term we have 
$${\E}\int_0^t\int_{|y|\leq 1}|y|\int_s^t g(u,s)(\frac{1-\beta}{u-s}+\lambda) \Big|\int_{[0,s)}\int_{|x|\leq 1}[\phi(u-v)-\phi(s-v)]x{\tilde N}(dv,dx) \Big|du\nu(dy)ds$$
$$\leq c\int_0^t\int_s^t g(u,s)(\frac{1-\beta}{u-s}+\lambda)\int_0^s |\phi(u-v)-\phi(s-v)|dvduds$$
$$\leq c\int_0^t\int_s^t g(u,s)(\frac{1-\beta}{u-s}+\lambda)(u-s)s^{\gamma}duds$$
$$\leq c\int_0^t\int_s^t (u-s)^{\gamma+\beta-1}(\frac{1-\beta}{u-s}+\lambda)duds<\infty.$$
Note that this last term is the only term that requires $\beta+\gamma>1.$ For the other terms, $\beta, \gamma>0$ is enough. 
Note that if $\beta\leq \frac{1}{2}$, the first mapping is not in $L^2$. \\
The study of the integrability of the second mapping is immediate assuming also $\gamma+\beta>1: $
$${\E}\int_0^t\int_{|x|\leq 1} |x|\int_s^t g(u,s)\phi(u-s)(\frac{1-\beta}{u-s}+\lambda)du\nu(dx)ds$$
$$\leq c\int_0^t\int_s^t g(u,s)\phi(u-s)(\frac{1-\beta}{u-s}+\lambda)duds.$$
Hence we conclude that the integral $\int_0^t Y(s-) dX(s)$ is well defined.
\end{example}

\begin{remark}
Following similar computations we can provide another example with $J$ a pure jump additive process by taking e.g.

$$\nu(dt,dx) = c h(t) \vert x \vert^{-1-\alpha} dtdx, \quad \alpha\in (0,2),$$ 
with $h$ a positive deterministic function such that $\int_0^t h(s)ds < \infty$.
\end{remark}

\section*{Acknowledgments}

\hspace{0.5cm} This work has been developed under the project Stochastic in Environmental and Financial Economics (SEFE) at the Center for Advanced Study (CAS) at the Norwegian Academy of Science and Letters. The authors thank CAS for the support and the kind hospitality.

\end{document}